%% file: comb.tex
\documentclass[12pt]{article}

\usepackage{amsfonts}
\usepackage[centertags]{amsmath}
\usepackage{amssymb}
\usepackage{amsthm}
\usepackage{newlfont}
\usepackage{graphics}
\usepackage{graphicx}
\usepackage{color}
\usepackage[latin1]{inputenc}
\usepackage[T1]{fontenc}
\usepackage[all]{xy}
\usepackage{a4,latexsym}
\input{vorspann20.tex}

\begin{document}

\title{A comb of origami curves in the moduli space $M_3$ with 
three dimensional closure}

\date{}
\author{{\it Frank Herrlich
 \footnote{
 e-mail: {\sf herrlich@math.uni-karlsruhe.de}, 
    }, \ 
    Gabriela Schmith\"usen}
  \footnote{
   e-mail: {\sf schmithuesen@math.uni-karlsruhe.de}}\\[3mm]
     \small
     Institut f\"ur Algebra und Geometrie, 
       Universit\"at Karlsruhe, 76128 Karlsruhe, Germany}

\maketitle
\begin{abstract}
\noindent
The first part of this paper is a survey on Teichm\"uller curves and Veech 
groups, with emphasis on the special case of origamis where much stronger
tools for the investigation are available than in the general case.\\
In the second part we study a particular configuration of origami curves
in genus 3: A ``base'' curve is intersected by infinitely many ``transversal''
curves. We determine their Veech groups and the closure of their locus 
in $M_3$, which turns out to be a three dimensional variety and the image 
of a certain Hurwitz space in $M_3$.
\end{abstract}

By an {\it origami} we mean a closed surface endowed with a translation 
structure that comes from glueing a certain number of squares. Variation 
of the 
translation structure leads to a 1-parameter family of Riemann surfaces, 
which in turn determines an algebraic curve in the moduli space of Riemann
surfaces. For every genus $g$, there are countably many such origami 
curves in $M_g$. So far, not much is known about the configurations of 
origami curves in general.\\[2mm]
In a previous paper \cite{WMS}, we found an origami $W$ of genus 3 for which
the origami curve is intersected by infinitely many other origami curves.
We call this configuration a ``comb'' although this name is misleading
in several respects: first, it may (and does) happen that a transversal
curve intersects the base curve in more than one point. And secondly, the
transversal curves do not all lie in one plane, not even on a surface:
in fact, one main result of this paper is that the closure of the set of
transversal origamis is a three dimensional algebraic subvariety of the
moduli space $M_3$ that we describe explicitly.\\[2mm]
As a second main result we determine the Veech groups of the
transversal origamis. We shall explain in Section \ref{vg} and \ref{veech}
how to define the Veech group of a translation surface and that, for an
origami $O$, it is a subgroup $\Gamma(O)$ of $\slzwei(\ZZ)$ of finite
index. Its importance
for the origami curve $C(O)$ lies in the fact that $C(O)$ is birational to
$\HH/\Gamma(O)$ (where $\HH$ is the upper half plane).\\[2mm]
The first two sections of this paper give a survey of Teichm\"uller curves
in general and origami curves as special cases. In particular we describe 
the Veech groups
corresponding to Teichm\"uller disks and to origamis. In Section \ref{wms}
we review the results of \cite{WMS} on the special origami curve $W$.
The last two sections contain the new results of this paper. In Section 
\ref{H} we determine the closure of the locus of the 
transversal origamis and relate it
to the Hurwitz space of certain coverings of elliptic curves. In Section
\ref{VG} we show that the Veech group of a transversal origami has index 3
in the stabilizer group of the corresponding configuration of critical
points. These are $n$-torsion points on an elliptic curve for varying $n$.
If $n$ is odd, we explicitly describe the Veech group.

\section{Teichm\"uller curves and Veech groups}
In this section we explain what {\em Teichm\"uller curves} are. They are 
special algebraic subvarieties of the moduli space $\mg$ of 
complex (regular) algebraic curves of genus $g$. They occur as 
images of one-dimensional analytic subvarieties of the  Teichm\"uller space
$\tg$, namely {\em Teichm\"uller disks}, 
which are natural with respect to the holomorphic structure on $\tg$
as well as to the Teichm\"uller metric.
Teichm\"uller curves are closely related to 
{\em Veech groups}, certain subgroups of $\slzwei(\RR)$. For a 
general introduction to Teichm\"uller spaces and the Teichm\"uller metric
see e.g. \cite{it}.

\begin{defn}
Let $\iota: \HH \into \tg$ be an embedding of the upper half
plane $\HH$ into Teichm\"uller space that is holomorphic and
isometric (with respect to the Poincar\'e metric on $\HH$ 
and the Teichm\"uller metric on $\tg$). 
\begin{enumerate}
\item[a)]
The image $\Delta$ of such an
embedding is called a {\em Teichm\"uller disk}.
\item[b)]
If $\Delta$ projects to an algebraic curve $C$ in the moduli space $\mg$, 
then this curve $C \subseteq \mg$ is called a {\em Teichm\"uller curve}.
\end{enumerate}
\end{defn}

First examples of Teichm\"uller curves were given by Veech 
in \cite{ve}. Since then, several
authors have been working on this subject. 
The reader can find a detailed overview and 
further hints to literature e.g. in \cite{hesh2}. In this section
we give a short introduction without proofs in order to provide the
reader with the general ideas and explain everything that one 
needs to understand the results  in Sections \ref{H} and \ref{VG}. For details
and proofs we refer the reader to \cite{hesh2}.\\

\subsection{Construction of Teichm\"uller disks}\label{td}
Teichm\"uller disks are obtained in the following way: Let $X$ be a Riemann
surface and $q$ a holomorphic quadratic differential. By integration $q$ 
naturally defines  a flat structure $\mu$ on $X\backslash\{\mbox{zeroes of } q\}$, 
i.e. an atlas such
that all transition maps are of the form
\[z \mapsto \pm z + c \;\; \mbox{ with } c \in \CC.\]
For each $A \in \slzwei(\RR)$, one gets a natural 
variation of the flat structure: 
Composing each chart with the affine map $z \mapsto Az$
defines a new flat structure $\mu_A$. 
We use here and in the whole article 
the $\RR$-linear identification of $\CC$ with $\RR^2$ that
maps $\{1,i\}$ to the standard basis of $\RR^2$.\\
Observe that the flat structure $\mu_A$ is in particular
a holomorphic structure on the topological surface underlying $X$ that will 
in general not be holomorphically equivalent to $\mu$. Furthermore,  the 
identity map defines a diffeomorphism $\id:X_{\mu} \to X_{\mu_A}$ 
between the Riemann surfaces 
$X_{\mu} := (X,\mu)$ and $X_{\mu_A} := (X,\mu_A)$ which is again in general
not holomorphic. Altogether one obtains 
a map
\[\iota: \slzwei(\RR) \to \tg, \quad  A \mapsto [X_{\mu_A},\id].\]
If $A$ is in $\sozwei(\RR)$, then $z \mapsto A\cdot z$ is in fact
biholomorphic, thus
the map $\iota$ factors through the quotient by $\sozwei(\RR)$ and one obtains the map
\[\overline{\iota}: \HH = \sozwei(\RR)\backslash\slzwei(\RR) \to \tg \quad;\]
one can show that $\bar{\iota}$ is a holomorphic and isometric embedding. Thus its image $\Delta = 
\Delta_{\mu}$ is a Teichm\"uller disk.

\subsection{Veech groups}\label{vg}

In order to study Teichm\"uller curves one needs to project 
the Teichm\"uller disk $\Delta = \Delta_{\mu}$ defined as above to the moduli space
$\mg$ by the natural projection $\proj:\tg \to \mg$. In general, the
Zariski closure of the image will be large, i.e. higher-dimensional. 
Only occasionally it is one-dimensional and $\proj(\Delta)$
is a Teichm\"uller curve $C$. The Veech group introduced in this subsection
makes it possible to decide whether this happens or not; and if it happens
the Veech group ``knows'' how the Teichm\"uller curve looks like.

\begin{defn}
Let $\Xstern$ be a Riemann surface together with a flat 
structure $\mu$.
\begin{itemize}
\item[a)]
A diffeomorphism $f: \Xstern\to \Xstern$ is called {\em affine}, if it is locally affine;
i.e. in terms of the charts of the flat atlas $\mu$ it is of the form:
\begin{equation}\label{affin}
z \mapsto Az + c \quad \mbox{ with }A \in \emglzwei(\RR) \mbox{ and } c \in \CC.
\end{equation}
\item[b)]
Observe that the matrix $A \in \emglzwei(\RR)$ is independent of the charts
up to the sign.
Let $\emaffplus(\Xstern,\mu)$ be the group of orientation preserving
affine diffeomorphisms and
\[D: \;\; \emaffplus(\Xstern,\mu) \to \empslzwei(\RR), \quad f \mapsto 
[A],\]
with $A$ as in a) and $[A]$ its class in $\empslzwei(\RR)$.
The {\em projective Veech group} $\Gammaquer(\Xstern,\mu)$ is the image 
of $D$ in $\empslzwei(\RR)$.
Hence, it contains the classes of all matrices that occur in {\em 
(\ref{affin})}
for orientation preserving affine diffeomorphisms $f$ on 
$(\Xstern,\mu)$.
\end{itemize}
\end{defn}

{\bf Remark:}\\[1mm]
{\bf 1.} For the examples of flat surfaces that we will study in this article, 
the flat structure is in fact a translation structure,
i.e. all transition maps are translations. 
In this case the matrices in 
(\ref{affin}) are in fact independent of the charts. With the notation 
from above one has a 
map $\affplus(\Xstern,\mu) \to \glzwei(\RR)$, which we also denote
by $D$. We call its image the 
{\em Veech group} $\Gamma(\Xstern,\mu)$. For $f \in 
\affplus(\Xstern,\mu)$,
we call $D(f) \in \Gamma(\Xstern,\mu)$ the {\em derivative} of $f$.\\
Furthermore, if $(\Xstern,\mu)$ is obtained from a closed 
Riemann surface $X$ together with a quadratic differential $q$ as described 
in \ref{td} with $\Xstern \subseteq X\backslash\{\mbox{zeroes of 
$q$}\}$, then 
$q$ defines a translation structure if and only if $q$ is the square of a 
holomorphic differential on $X$. In this case
$\Gamma(\Xstern,\mu)$ is a subgroup of 
$\slzwei(\RR)$. We 
denote the group of orientation preserving affine
diffeomorphisms of $X^*$ also by $\affplus(X,\mu)$, 
the Veech group 
by $\Gamma(X,\mu)$ and the projective Veech group by 
$\Gammaquer(X,\mu)$.\\[1mm] 
{\bf 2.} One may observe  from the 
definitions that \[\Gamma(X,\mu_A) = A\cdot\Gamma(X,\mu)\cdot A^{-1}, \] 
where the flat structure $\mu_A$ is defined as in Section \ref{td}. Hence,
Veech groups for flat surfaces on the same Teichm\"uller disk are
equal up to conjugation in $\slzwei(\RR)$.\\

Recall that the natural projection $\tg \to \mg$ is the quotient map by the
mapping class group
\[\Gamma_g := \diffplus(X)/\diffplus_0(X),\]
where $\diffplus(X)$ is the group of orientation preserving diffeomorphisms and
$\diffplus_0(X)$ consists of the diffeomorphisms homotopic to the identity.\\
Let $\Delta = \Delta_{\mu}$ be a Teichm\"uller disk as in Section \ref{td}.
For our purposes, the following well known properties of the affine group 
and the Veech group $\Gamma$ of $(X,\mu)$  are important:
\begin{itemize}
\item 
The affine diffeomorphisms are precisely the diffeomorphisms that
stabilize $\Delta$, i.e. map $\Delta$ to itself.
\item
The kernel of $D: \affplus(X,\mu) \to \Gamma \; \subseteq  \;\slzwei(\RR)$ consists
of those diffeomorphisms that act trivially on $\Delta$.
Consequently, one has:
\begin{equation}\label{stabs}
\Gamma \;\; \cong \;\; \affplus(X,\mu)/\mbox{kernel of D} \;\; \cong \;\;
\mbox{global stabilizer}/\mbox{pointwise stabilizer} 
\end{equation}
\item By the identification $\iota: \HH \to \Delta$, the action of $\Gamma$ on
$\Delta$ given by (\ref{stabs}) becomes equal to the action of 
$\Gammastern :=  R\Gamma R^{-1} \subseteq \slzwei(\RR)$ on $\HH$ as Fuchsian group with
\[ R = \bpm-1&0\\0&1\epm.\]
\end{itemize}
Thus the map $\proj: \Delta \to \proj(\Delta) \subseteq \tg$ factors through
 $\HH/\Gammastern$.\\
Observe that $\HH/\Gammastern$ is an algebraic curve if and only if the quotient has finite
volume, i.e. if $\Gamma$ is a lattice in $\slzwei(\RR)$. We put these results
together in:\\[3mm]
{\bf Theorem A\footnote{For a general introduction to algebraic curves and algebraic geometry 
we refer to \cite{ha}.}}
{\it 
Let $\Delta$ be a Teichm\"uller disk in the Teichm\"uller space $\tg$.
\begin{itemize}
\item
The image $\emproj(\Delta)$ is an algebraic curve $C$ in $\mg$, if and only if 
the Veech group $\Gamma$ is a lattice in $\emslzwei(\RR)$.
\item In this case, $\HH/\Gammastern$ is the normalization of the 
algebraic curve $C$.
\end{itemize}}
Hence, the Veech group detects, whether a Teichm\"uller disk leads to a 
Teichm\"uller curve and, if this happens, it 
determines the Teichm\"uller curve up to birationality.

\section{Origamis}

One way to obtain particular flat surfaces is given by the following construction:\\
Take finitely many unit squares and glue their edges by translations, such that:
\begin{itemize}
\item
Each left edge is glued to a right one and vice versa.
\item
Each upper edge is glued to a lower one and vice versa.
\item 
The resulting surface $X$ is connected.
\end{itemize}
$X$ carries in a natural way a flat structure defined by the unit squares. 
Actually one has even more, namely a translation structure, 
since the glueings are done by translations.\\
Surfaces obtained this way are called {\em square tiled surfaces}
or {\em origamis} as Pierre Lochak baptized them in \cite{l} due to the
playful character of their definition.\\
The baby example is to take only one square. There is only one 
possibility to glue its edges according to the rules and one obtains a torus.
In fact this example is universal in the following sense:
Each other surface $X$ obtained by the construction above
allows a morphism to the torus by mapping
each square on $X$  to the unique square of the torus. 
This morphism is unramified except for the one point on
the torus that comes from the vertices of the square. Conversely,
 a surface that allows such a 
covering can be obtained by glueing squares according to the rules above.  This 
motivates the following definition
for (oriented) origamis. 
\begin{defn}\label{deforigami}
An {\em origami} $O$ is a covering $p:X \to E$ from a topological surface $X$ 
to the torus $E$ that is ramified at most over one point $\infty \in E$.
 Two origamis
$O = (p:X \to E)$ and $O' = (p':X' \to E)$ are said to be {\em isomorphic},
if there is some homeomorphism $f:X \to X'$ such that $p' \circ f = p$.
\end{defn}
One then obtains the flat structure $\mu$ on $X$ equivalently as follows:
\begin{itemize}
\item Identify the torus $E$ with $\CC/\Lambda_I$, where
$\Lambda_I$ is the lattice in $\CC$ generated by $1$ and $i$.
\item Take the natural translation structure on $E = \CC/\Lambda_I$ coming from 
locally reversing the universal covering $\CC \to E$.
\item
Lift it by the unramified covering $p: X^* = X-p^{-1}(\infty) \to E^* = E-\{\infty\}$
to  the translation structure $\mu$ on $X^*$.
\end{itemize}
Obviously, isomorphic origamis define isomorphic translation surfaces.
Note that for origamis the translation structure comes from
the holomorphic quadratic differential $\omega^2 = p^*(\omega_E^2)$ on $X$, 
where
$\omega_E$ is the (essentially unique) holomorphic differential on $E$.

\subsection{Description of origamis}
\label{descr}

One advantage of flat surfaces coming from origamis is that they are 
described by only few combinatorial data: the {\em glueing data} 
that record which square is glued to which one (in horizontal resp.\ in
vertical direction). The equivalent description in terms of Definition 
\ref{deforigami} is the monodromy 
\[\pi_1(E^*) \to S_d,\] 
where $d$ is the degree
of $p$ and $S_d$ is the symmetric group. Observe that 
the fundamental group $\pi_1(E^*)$ is the free group in two  generators. 
The images 
in $S_d$ of the two generators $x$, the horizontal closed curve on $E^*$,
and y, the vertical closed geodesic, give precisely the glueing data.\\
For a third equivalent description pointed out in \cite{sh1} observe that,
by the theorem of the universal covering, 
the unramified covering $X^* \to E^*$ gives an embedding
\begin{equation}\label{uinf2}
U = \pi_1(X^*) \;\; \into \;\;\pi_1(E^*) = F_2 
\end{equation}
of index $d$ equal to the $\mbox{degree of } p$. Vice versa one obtains, 
again by the theorem of the universal covering, to each such subgroup
$U$ of $F_2$ a covering $p$ as in Definition \ref{deforigami}. Hence, origamis
correspond to finite index subgroups of $F_2$.

\subsection{Veech groups of origamis}
\label{veech}

A second advantage of origamis is that one has an alternative approach
to their Veech groups, which we describe in this subsection. 
It makes them easier to handle than those of general 
flat surfaces.\\ 

Let $O$ be an origami and $(X,\mu)$ the flat surface that it defines.
\begin{defn}
The Veech group of the origami $O$ is the Veech group of $(X,\mu)$
\[\Gamma(O) := \Gamma(X,\mu)\]
\end{defn}
One may start with the baby example from above: 
the origami with only one square; 
in terms of Definition~\ref{deforigami} it is the covering $O = (\id:E \to E)$. 
The flat surface then is $\CC/\Lambda_I$. Its affine diffeomorphisms come from 
affine diffeomorphisms of $\CC$ that descend to $\CC/\Lambda_I$. Hence, they are 
precisely the affine diffeomorphisms preserving the lattice $\Lambda_I = \ZZ \oplus \ZZ i$ .
Thus $\Gamma(O) = \slzwei(\ZZ)$.\\
For a general origami $O$ the Veech group is always a finite index subgroup of $\slzwei(\ZZ)$:
\[\Gamma(O) \subseteq \slzwei(\ZZ)\]
The converse also holds: By the Theorem of Gutkin and Judge in \cite{gj}
a Veech group is a finite index subgroup of $\slzwei(\ZZ)$ if and only if the translation 
surface
comes from a once ramified covering of a torus.\\

One may now use the above description of origamis by finite index subgroups of 
$F_2$ 
 to determine the Veech group. In \cite{sh1}, the following
characterization was given:\\

{\bf Theorem B:}
{\it
\begin{itemize}
\item 
Let $U$ be the finite index subgroup of $F_2$ defined in (\ref{uinf2}) for the origami $O$.
\item 
Let $\betadach: \emaut(F_2) \to \emglzwei(\ZZ) = \emout(F_2)$ be the natural projection and \\
$\emautplus(F_2) := \betadach^{-1}(\emslzwei(\ZZ))$.
\item
Let $\emstab(U) := \{\gamma \in \emautplus(F_2)| \gamma(U) = U\}$ be the 
{\em stabilizing group of $U$}
\end{itemize} 
Then for the Veech group $\Gamma(O)$ holds:  $\Gamma(O) = \betadach(\emstab(U))$.
}\\

{\bf Remark:}
The proof of Theorem B is based on the theorem of the universal covering.
One uses a fixed universal covering 
$u:\HH \to E^*_{I} = \CC/\Lambda_I - \{\infty\}$. 
One then considers the group 
\[\aff(\HH,u) := \{\tilde{f} \in \diffplus(\HH)|\, \tilde{f} 
\mbox{ descends via $u$ 
to an affine map on $E^{*}_I$}\}.\] In fact, this is the group of affine
diffeomorphisms of $\HH$ with respect to the translation structure that comes 
as lift via $u$  from the one on  $E^*_{I}$. 
One then has the isomorphism:
\[*: \aff(\HH,u) \to \autplus(F_2), \quad 
      \tilde{f} \mapsto (\gamma: w \mapsto \tilde{f}w\tilde{f}^{-1}). \]
Here we use that $F_2 = \gal(\HH/E^*_I)$ is the
group of deck transformations.
If one now has  an origami $O = (p:X\to E)$, one may restrict to those affine 
diffeomorphisms on $\HH$ that descend to $X^*$ and obtains the following commutative
diagram:
\[ \xymatrix{
   \{ \tilde{f} \in \aff(\HH,u)| \tilde{f} \mbox{ descends to } \Xstern\} 
       \ar[d]_{D} \ar[rr]^{\hspace*{20mm}\star}  &&
      \stab(U) \ar[d]^{\betadach}&\\
      \Gamma(O)  \ar[rr]^{\scs \id} && \betadach(\stab(U)) &\subseteq \slzwei(\ZZ) 
    }
\]
The map $\betadach:\autplus(F_2) \to \slzwei(\ZZ)$ can be given explicitly as follows:
Let  $\sharp_x(w)$ be the number of occurrences of $x$ in $w$ (with $x^{-1}$ counted 
as $-1$) and 
similarly for $\sharp_y(w)$. Then
\[\betadach(\gamma) =  \bpm \sharp_x(\gamma(x)) & \sharp_x(\gamma(y))\\ 
                      \sharp_y(\gamma(x)) & \sharp_y(\gamma(y)) \epm.\]

\subsection{Teichm\"uller curves of origamis}
\label{ocurves}

A further advantage of flat surfaces coming from origamis is 
that they always define Teichm\"uller curves: As we have mentioned in the last
subsection, their Veech groups are finite index subgroups of $\slzwei(\ZZ)$
and thus lattices in $\slzwei(\RR)$. Hence the Teichm\"uller disk $\Delta$ defined by
such a surface always projects to a Teichm\"uller curve in the moduli space (see 
Section~\ref{vg}).\\

Observe that composing each chart of the translation surface $\CC/\Lambda_I$
with the affine map $z \mapsto Az$ with $A \in \slzwei(\RR)$ is equivalent 
to taking the translation surface $\CC/\Lambda_A$, where $\Lambda_A$ is the lattice 
\[\Lambda_A = <\bpm a\\c\epm , \bpm b\\d \epm> \quad \mbox{ for } A = \bpm a&b\\c&d\epm.\]
Here again we identify $\CC$ with $\RR^2$.\\
Therefore one may describe the Teichm\"uller disk $\Delta$ 
defined by  an origami $O = (p:X\to E)$ as follows:
\begin{itemize}
\item 
For each $A \in \slzwei(\RR)$ let $\eta_A$ be the translation structure on $E$ obtained by identifying
$E$ with $\CC/\Lambda_A$.
\item
Let $\mu_A$ be the lift of $\eta_A$ to $X^* = X\backslash p^{-1}(\infty)$
and call $X_A = (X,\mu_A)$.
\end{itemize}
Then the Teichm\"uller disk $\Delta$ is given as
\[\Delta =\{[X_A,\id]| \; A \in \slzwei(\RR)\},\]
where, as before, $\id$ is  topologically the identity on $X$, 
whereas the holomorphic structure on $X$ is changed.\\

We denote the Teichm\"uller curve defined by an origami $O = (p:X\to E)$ 
by $C(O)$ and call it  an {\it origami curve}; recall that
it is the image in the moduli space of the Teichm\"uller disk $\Delta$.
It is possible to decide whether two origamis induce the same origami curve 
in terms of the subgroup $U\subset F_2$ associated to an origami in Section
\ref{descr}:
\begin{prop}
\label{gleich}
Let $O=(p:X\to E)$ and $O'=(p':X'\to E)$ be two origamis and $U$ 
(resp.\ $U'$) the corresponding subgroups of $F_2$. Then\\[1mm]
{\bf a)} $O$ is isomorphic to $O'$ if and only if $U$ is conjugate to $U'$
in $F_2$.\\[1mm]
{\bf b)} $C(O)=C(O')$ if and only if there is a $\gamma\in\mbox{\em Aut}^+(F_2)$ 
such that $\gamma(U)=U'$.
\end{prop}
\begin{proof}
a) By definition, $O$ is isomorphic to $O'$ if and only if there is a
diffeomorphism $f:X\to X'$ satisfying $p=p'\circ f$. Then also the unramified
coverings $X^*\to E^*$ and $(X')^*\to E^*$ are isomorphic, and hence
$\pi_1(X^*)$ is conjugate to $\pi_1(X')^*$ in $\pi_1(E^*)$.\\[1mm]
b) $C(O)$ is equal to $C(O')$ if and only if there is a
diffeomorphism $f:X\to X'$ such that $p^*(\omega_E)$ is equal to 
$f^*((p')^*(\omega_E))$ up to a multiplicative constant in $\CC$, where
$\omega_E$, as before, is the holomorphic differential on $E = E_I$.
Equivalent conditions are that $f$ descends to an affine diffeomorphism of $E$
or that $f$ lifts to an element $\tilde{f} \in \aff(\HH,u)$ in the notation
of the remark following Theorem~B. Under the isomorphism $*$ in the same remark, 
such an $f$ 
corresponds to an
automorphism of $F_2$ that maps $U$ to $U'$.
\end{proof}

\input{wms.tex}

\input{vg.tex}

\end{document}

{\color{blue} }

%% file: vorspann20.tex
\newtheorem{thm}{Theorem} 
\newtheorem{cor}{Corollary} 
 
\newtheorem{prop}[cor]{Proposition} 
\newtheorem{lemma}[cor]{Lemma}
\newtheorem{defn}[cor]{Definition}

\newcounter{diagramm}

\newcommand{\fs}{\footnotesize}
\newcommand{\scs}{\scriptsize}

\newcommand{\CC}{\mathbb{C}}

\newcommand{\RR}{\mathbb{R}}
\newcommand{\ZZ}{\mathbb{Z}}
\newcommand{\NN}{\mathbb{N}}
\newcommand{\HH}{\mathbb{H}} 

\newcommand{\FF}{\mathbb{F}}
\newcommand{\AAA}{\mathbb{A}}
\newcommand{\PP}{\mathbb{P}}

\newcommand{\gal}{\mbox{Gal}}

\newcommand{\emaut}{\mbox{\em Aut}}  
\newcommand{\autplus}{\mbox{Aut}^+}
\newcommand{\emautplus}{\mbox{\em Aut}^+}

\newcommand{\id}{\mbox{id}}

\newcommand{\aff}{\mbox{Aff}}
\newcommand{\affplus}{\mbox{Aff}^+}
\newcommand{\emaffplus}{\mbox{\em Aff}^+}
\newcommand{\emaff}{\mbox{\em Aff}}

\newcommand{\slzwei}{\mbox{SL}_2}
\newcommand{\emslzwei}{\mbox{\em SL}_2}

\newcommand{\empslzwei}{\mbox{\em PSL}_2}
\newcommand{\glzwei}{\mbox{GL}_2}

\newcommand{\emglzwei}{\mbox{\em GL}_2} 

\newcommand{\sozwei}{\mbox{SO}_2}

\newcommand{\emout}{\mbox{\em Out}}

\newcommand{\proj}{\mbox{proj}}
\newcommand{\emproj}{\mbox{\em  proj}}

\newcommand{\mg}{M_g}
\newcommand{\tg}{T_g}
\newcommand{\diff}{\mbox{Diff}}
\newcommand{\diffplus}{\diff^+}

\newcommand{\into}{\hookrightarrow}

\newcommand{\bpm}{\begin{pmatrix}}
\newcommand{\epm}{\end{pmatrix}}
\newcommand{\Gammastern}{\Gamma^*}
\newcommand{\Gammaquer}{\bar{\Gamma}}

\newcommand{\stab}{\mbox{Stab}}
\newcommand{\emstab}{\mbox{\em Stab}}
\newcommand{\betadach}{\hat{\beta}}

\newcommand{\el}{E_{-1}}
\newcommand{\trnstern}{E^{[n]*}}
\newcommand{\trstern}{E^{*}}
\newcommand{\Xstern}{X^*}

\newcommand{\ifff}{\Leftrightarrow}

%% file: wms.tex
\newcommand{\Aut}{\mbox{Aut}}
\newcommand{\emAut}{\mbox{\em Aut}}

\section{The quaternion origami}
\label{wms}
In this section we illustrate the general concepts that we introduced in the 
first two sections by a specific example of an origami curve in genus 3. 
In addition to the general features shared by all origami curves, our origami 
$W$ has several remarkable properties which make it particularly 
interesting. We restrict ourselves to a short description and refer to
\cite{WMS} for more details and, in particular, for proofs.\\[2mm]
\subsection{$W$ as an origami}
\label{kaestchen}
Using the combinatorial definition of origamis, $W$ can be described by 8 
squares that are glued as indicated (edges are glued if they have the same 
label).

\begin{center}
\setlength{\unitlength}{.8cm}
\begin{minipage}{13cm}
\begin{picture}(5,3)
\hspace*{-4cm}
\put(11,0){\framebox(1,1)}
\put(13,0){\framebox(1,1)}
\put(11,1){\framebox(1,1){1}}
\put(12,1){\framebox(1,1)}
\put(13,1){\framebox(1,1)}
\put(14,1){\framebox(1,1)}
\put(12,2){\framebox(1,1)}
\put(14,2){\framebox(1,1)}

\put(11.3,-.1){\scriptsize{//}}
\put(13.4,-.1){\scriptsize{/}}
\put(12.3,.9){\scriptsize{$\backslash\backslash$}}
\put(14.4,.9){\scriptsize{$\backslash$}}
\put(11.4,1.9){\scriptsize{/}}
\put(13.3,1.9){\scriptsize{//}}
\put(12.4,2.9){\scriptsize{$\backslash$}}
\put(14.3,2.9){\scriptsize{$\backslash\backslash$}}

\put(10.83,.28){=}
\put(10.83,.45){=}
\put(11.9,.4){--}
\put(12.83,.4){=}
\put(13.83,.45){=}
\put(13.9,.33){--}
\put(11.9,2.4){--}
\put(12.83,2.4){=}
\put(13.83,2.45){=}
\put(13.9,2.33){--}
\put(14.83,2.28){=}
\put(14.83,2.45){=}
\put(10.87,1.35){$\circ$}
\put(14.87,1.35){$\circ$}

\put(10.9,-.1){$\bullet$}
\put(12.9,-.1){$\bullet$}
\put(10.9,1.9){$\bullet$}
\put(12.9,1.9){$\bullet$}
\put(14.9,1.9){$\bullet$}
\put(11.9,-.1){$\ast$}
\put(13.9,-.1){$\ast$}
\put(11.9,1.9){$\ast$}
\put(13.9,1.9){$\ast$}
\put(10.9,.9){\framebox(.2,.2)}
\put(12.9,.9){\framebox(.2,.2)}
\put(14.9,.9){\framebox(.2,.2)}
\put(12.9,2.9){\framebox(.2,.2)}
\put(14.9,2.9){\framebox(.2,.2)}
\put(11.8,.8){$\Diamond$}
\put(11.8,2.8){$\Diamond$}
\put(13.8,.8){$\Diamond$}
\put(13.8,2.8){$\Diamond$}

\end{picture}
\end{minipage}
\end{center}

Note that every vertex is glued to precisely one of the vertices of the square 
labeled ``1''. Euler's formula gives $2 - 2g = 8 - 16 + 4$, thus the genus 
of $W$ is 3. \\[2mm]
The total angle at every vertex is $4\pi$, so they are all ramification points 
of order 2 for the covering $p:W\to E$ of degree 8 to the torus $E$. 
Recall that $p$ is
obtained by mapping each of the eight squares to $E$. The map $p$
can be decomposed as follows: Observe that ``translation by 2 to the right''
is an automorphism of $W$, and let $q:W\to \overline{W}$ be the quotient map for
this automorphism. Then $\overline{W}$ is the origami\\[3mm]
\setlength{\unitlength}{.8cm}
\hspace*{2cm}
\begin{minipage}{5cm}
\begin{picture}(5,3.5)
\put(3,0){\framebox(1,1)}
\put(3,1){\framebox(1,1)}
\put(4,1){\framebox(1,1)}
\put(4,2){\framebox(1,1)}

\put(6.5,1.6){which is}
\put(6.5,1){equivalent to}

\put(10,.5){\framebox(1,1)}
\put(11,.5){\framebox(1,1)}
\put(10,1.5){\framebox(1,1)}
\put(11,1.5){\framebox(1,1)}

\put(2.9,.4){=}
\put(3.9,.4){--}
\put(3.9,2.4){--}
\put(4.9,2.4){=}
\put(3.3,-.1){\scriptsize{//}}
\put(3.3,1.9){\scriptsize{//}}
\put(4.4,.9){\scriptsize{/}}
\put(4.4,2.9){\scriptsize{/}}
\put(2.87,1.35){$\circ$}
\put(4.87,1.35){$\circ$}

\put(9.9,1.9){=}
\put(9.9,.9){--}
\put(11.9,.9){--}
\put(11.9,1.9){=}
\put(11.4,.4){\scriptsize{//}}
\put(11.4,2.4){\scriptsize{//}}
\put(10.5,.4){\scriptsize{/}}
\put(10.5,2.4){\scriptsize{/}}

\end{picture}
\end{minipage}
\begin{center}
\end{center}

Thus $\overline{W}$ is the torus $E$. 
{ As in Section 2, we identify $E$ with $\CC/\Lambda_I$
in such a way that the vertex of the square becomes the origin of the
group structure. In this identification, the 4 vertices of $W$ are mapped 
by $q$ to the points of order $2$ on $E$, and $p$ is obtained by postcomposing 
$q$ with the multiplication by 2 on $E$. As usual (see e.\,g.\ \cite{sil})
we call, for $n\ge2$, the points of order dividing $n$ in the group
$\CC/\Lambda_I$ the {\it $n$-torsion points} of $E$. They form a group 
$E[n]$ of order $n^2$, isomorphic to $\ZZ/n\ZZ\times\ZZ/n\ZZ$; in other words,
$E[n]$ is the kernel of the multiplication $[n]$ by $n$ on $E$.}

\subsection{The automorphism group of $W$}
\label{autos}
$W$ is a Galois origami, i.\,e.\ the covering $p:W\to E$ is a normal covering.
The Galois group is the quaternion group $Q$ of order 8. 
This can be seen e.\,g.\ by labeling the squares of $W$ as follows:

\begin{center}
\setlength{\unitlength}{.8cm}
\begin{minipage}{13cm}
\begin{picture}(5,3)
\hspace*{-4cm}
\put(11,0){\framebox(1,1){$-j$}}
\put(13,0){\framebox(1,1){$j$}}
\put(11,1){\framebox(1,1){1}}
\put(12,1){\framebox(1,1){$i$}}
\put(13,1){\framebox(1,1){$-1$}}
\put(14,1){\framebox(1,1){$-i$}}
\put(12,2){\framebox(1,1){$k$}}
\put(14,2){\framebox(1,1){$-k$}}

\put(11.3,-.1){\scriptsize{//}}
\put(13.4,-.1){\scriptsize{/}}
\put(12.3,.9){\scriptsize{$\backslash\backslash$}}
\put(14.4,.9){\scriptsize{$\backslash$}}
\put(11.4,1.9){\scriptsize{/}}
\put(13.3,1.9){\scriptsize{//}}
\put(12.4,2.9){\scriptsize{$\backslash$}}
\put(14.3,2.9){\scriptsize{$\backslash\backslash$}}

\put(10.83,.28){=}
\put(10.83,.45){=}
\put(11.9,.4){--}
\put(12.83,.4){=}
\put(13.83,.45){=}
\put(13.9,.33){--}
\put(11.9,2.4){--}
\put(12.83,2.4){=}
\put(13.83,2.45){=}
\put(13.9,2.33){--}
\put(14.83,2.28){=}
\put(14.83,2.45){=}
\put(10.87,1.35){$\circ$}
\put(14.87,1.35){$\circ$}

\end{picture}
\begin{center}
\end{center}
\end{minipage}
\end{center}

As usual, the elements of $Q$ are denoted $\pm1,\pm i,\pm j$ and $\pm k$,
with relations $i^2 = j^2 = k^2 = -1$ and $ij = k = -ji$.\\[2mm]
The full automorphism group $G$ of $W$ consists of the affine diffeomorphisms
with derivative $I$ and $-I$. The first ones are precisely the elements of $Q$.
$G$ is an extension of $Q$ of degree 2; more precisely, $G$ is 
a group of order 16 that contains (besides the
elements of $Q$) 6 involutions which act on $W$ by rotating the squares by an 
angle of $\pi$; each of these rotations has four fixed points which are 
midpoints either of squares or of edges. The last two elements of $G$, $c$
and $c^{-1}$, are the rotation by $\pi$ around the four vertices, and its 
inverse; note that $c$ has order 4 since the total angle at a vertex is $4\pi$.
It is a useful observation that $c^2$ equals the translation by 2 squares
(in horizontal or vertical direction: they are equal!); thus $c^2$ is the
element $-1$ in $Q$. Moreover, $c$ generates the center of $G$.

\subsection{The equation of $C(W)$}
\label{equ}
In this section we determine the origami curve $C(W)$ by explicitly giving
the algebraic equation of the covering Riemann surfaces.
Since $c$ has 4 fixed points, it follows from the Riemann-Hurwitz formula 
that $W/<c>$ has genus 0. Thus $W\to W/\!<\!c\!>$ is a cyclic covering of 
degree 4 of
the projective line, branched over 4 points which we can normalize to
be $0$, $1$, $\infty$ and some parameter 
$\lambda \in \PP^1 - \{0,1,\infty\}$. Therefore, as an algebraic curve, 
$W$ is given by an equation of the form
\begin{equation}
\label{lambda}
y^4 = x^{\varepsilon_0}(x-1)^{\varepsilon_1}(x-\lambda)^{\varepsilon_\lambda}.
\end{equation}
The exponents $\varepsilon_0$, $\varepsilon_1$ and $\varepsilon_\lambda$
are determined by the monodromy (which can be read off from the origami) 
and turn out to be 1. The complex structure on $W$ is determined by
$\lambda$ up to replacing $\lambda$ by one of $\frac{1}{\lambda}$, $1-\lambda$,
$1-\frac{1}{\lambda}$, $\frac{\lambda}{\lambda-1}$ and $\frac{1}{1-\lambda}$.
The covering $q:W\to E$ from Section \ref{kaestchen} is given by 
$(x,y)\mapsto(x,y^2)$, and $E$
is the elliptic curve with equation $y^2=x(x-1)(x-\lambda)$.\\[1mm]
{ In the last equation, $\lambda$ is the {\it Legendre parameter}
of the elliptic curve $E$. It is related to the lattice description 
$E = \CC/\Lambda_A$, $A\in\slzwei(\RR)$, in Section 2.3 by the classical 
modular function $\lambda$ on the upper half plane, which is the universal 
covering $\HH \to \PP^1 - \{0,1,\infty\}$ (see e.\,g.\ \cite[Chap.~VI.3]{FL}):
For $A = \begin{pmatrix}a&b\\c&d\end{pmatrix}\in\slzwei(\RR)$, the complex
torus $\CC/\Lambda_A$ is biholomorphic to the elliptic curve $E_{\lambda_A}$
with affine equation $y^2 = x(x-1)(x-\lambda_A)$, where 
$\lambda_A = \lambda(\frac{d\cdot i+b}{c\cdot i+a})$. For example, the unit 
lattice $\Lambda_I$ corresponds to $\lambda_I=-1$ (or to $\lambda=2$ or
$\lambda=\frac{1}{2}$, if we choose a different basis for $\Lambda_I$).\\[2mm]
The holomorphic differential $\omega = p^*(\omega_E)$ that according to
the introduction of Section 2 determines the translation structure on $W$,
can also be written down explicitly: it is}
{
\begin{equation}
\omega = p^*(\omega_E) = q^*[2]^*(\frac{dx}{2y}) = q^*(\frac{dx}{y}) = 
\frac{dx}{y^2}.
\end{equation}
}
\\[2mm]
An equivalent form of the equation (\ref{lambda}) is
\begin{equation}
\label{a}
y^4 = x^4 +2ax^2 +1
\end{equation}
where now $a$ is in $\PP^1 - \{1,-1,\infty\}$. The two equations are 
transformed into each other by a projective change of coordinates, see e.\,g.\
\cite{gu}.
The relation between $\lambda$ and $a$ is given by
\begin{equation}
a = \frac{\lambda + 1}{\lambda - 1},\qquad \lambda = \frac{a+1}{a-1}.
\end{equation}
The result of Section \ref{autos} states that, for every admissible value
of $a$ (or $\lambda$), the projective curve given by (\ref{a}) (resp.\ by 
(\ref{lambda})) has an automorphism group containing $G$. There are precisely 
two curves in the family that have a larger automorphism group: one is the
Fermat curve $y^4 = x^4 + 1$, which has 96 automorphisms. The other 
exceptional curve occurs for $\lambda = \frac{1}{2}\pm\frac{i}{2}\sqrt{3}$;
it has 48 automorphisms. { See \cite{KK} for a detailed 
discussion of these automorphism groups.} 

\subsection{The Veech group of $W$}
\label{veechwms}
We have seen in Section \ref{veech} that the Veech group of any origami is
a subgroup of $\slzwei(\ZZ)$ of finite index. For the origami $W$ we have:
\begin{prop}
$\Gamma(W) = \emslzwei(\ZZ)$.
\end{prop}
This can be checked by showing that the two generators of $\slzwei(\ZZ)$ -- 
e.\,g.\ the matrices $T = \begin{pmatrix}1&1\\0&1\end{pmatrix}$ and 
$S = \begin{pmatrix}0&-1\\1&0\end{pmatrix}$ --
are in fact affine automorphisms of $W$. A more conceptual argument is the 
following:\\
Let $U = \pi_1(W^*)$ be the subgroup of $F_2$ that corresponds to the
origami $W$ as in Section \ref{descr}.
By Section \ref{autos}, $U$ is the 
kernel of a surjective homomorphism $F_2\to Q$. It is not hard to see that all 
surjective homomorphisms $F_2\to Q$ have the same kernel. Therefore $U$ is
invariant under all automorphisms of $F_2$, which by Theorem B
implies that $\Gamma(W) = \slzwei(\ZZ)$.\\[2mm]
As a consequence, $\HH/\Gamma(W)$ is isomorphic to the affine line. It turns 
out that, in this case, the map $\HH/\Gamma(W)\to C(W)$ to the origami curve in
$M_3$ (which by Theorem~A is birational in general) is an isomorphism. 
Thus we have
\begin{prop}
The origami curve $C(W)$ is an affine line, embedded into the moduli space 
$M_3$.
\end{prop}

\subsection{Origami curves crossing $W$}
\label{cross}
The origami $W$ has several other remarkable properties. One of them is that 
the Jacobian of each of the Riemann surfaces $W_\lambda$, 
$\lambda\in\PP^1-\{0,1,\infty\}$, splits (up to isogeny) into a product of 
three elliptic curves, see \cite[Prop.~1.6]{WMS}. More precisely, 
$\mbox{Jac}(W_\lambda)$ is isogenous to $E_\lambda\times E_{-1}\times E_{-1}$,
where $E_\lambda$ is the elliptic curve with equation $y^2 = x(x-1)(x-\lambda)$
and $E_{-1}$ (given by the equation $y^2 = x^3-x$) is independent of $\lambda$;
{ see \cite[Chap.~XI]{BL} for background information on Jacobian
varieties of Riemann surfaces.}
The fact that the Jacobian of the family of curves defined by $W$ has a
constant part of codimension 1 shows that $C(W)$ is a Shimura curve,
{ see \cite{M1}. In that paper, M\"oller also shows} that 
among all Teichm\"uller curves that are induced 
by the square of a holomorphic differential, $C(W)$ is the only one that is
at the same time a Shimura curve.\\[2mm]
For the present note, the following result on $C(W)$, proved in 
\cite[Thm.~3.1]{WMS}, is the most interesting:\\[2mm]
{\bf Theorem C:} 
{\it $C(W)$ intersects infinitely many other origami curves.}\\[2mm]
Since all origami curves are defined over number fields, there can be at most
countably many origami curves that intersect a given one. In particular, their
image in the moduli space $M_3$ of Riemann surfaces of genus 3 cannot be
closed. We shall study the closure of the set of origami curves intersecting 
$C(W)$ in the next section. Here we briefly recall how these origamis 
arise:\\[2mm]
We have seen in Section \ref{autos} that, for all 
$\lambda\in\PP^1-\{0,1,\infty\}$, $\Aut(W_\lambda)$ contains 7 involutions.
All of them have 4 fixed points on $W_\lambda$, so the quotient curve has
genus 1 in all cases. One involution is the translation by 2 squares which
induces the map $q$ to $E_\lambda$, see Section \ref{kaestchen}. For all
the other involutions the central automorphism $c\in\Aut(W_\lambda)$ 
descends to an automorphism $\bar c$ of order 4 on the quotient. Since
$\bar c$ has two fixed points, which we call $O$ and $M$,
the quotient must be the elliptic curve 
$E_{-1}$, and $\bar c$ is a rotation around $O$ and $M$ of angle
$\frac{\pi}{2}$.
\\[2mm]
Let us denote by $\sigma$ one of these six involutions and by 
$\kappa = \kappa_\lambda:W_\lambda\to E_{-1}$ the quotient map. 
The ramification points of $\kappa$ are the four fixed points of $\sigma$. 
Their images
$P_0(\lambda),\dots,P_3(\lambda)$ form an orbit under $\bar c$. Choose the
fixed point $O$ of $\bar c$ as
the origin on $E_{-1}$. 
Then if $P_0(\lambda)$ is an $n$-torsion point
for some $n$, the same holds for $P_1(\lambda)$, $P_2(\lambda)$ and
$P_3(\lambda)$. Hence in this case, 
$[n]\circ\kappa$ is ramified only over the origin 
and thus defines an origami. We showed in \cite{WMS}
by explicit calculations that, for { any} $n\ge 3$ and 
{ any} $n$-torsion point
$P$ on $E_{-1}$,
there is a $\lambda\in\PP^1-\{0,1,\infty\}$ such that
$P_0(\lambda)=P$.

\section{The Hurwitz space $H$}
\label{H}

\subsection{The origamis $D_P$}
\label{dp}
We have seen that, for each $n\ge3$, every $n$-torsion point $P$ on the
elliptic curve $E_{-1}: y^2 = x^3-x$ induces an origami $D_P$ of degree $2n^2$
whose origami curve $C(D_P)$ intersects the origami curve $C(W)$.
Combinatorially, $D_P$ consists of two large squares made of $n^2$ small 
squares each; the point $P =P_0$ corresponds to a (primitive) vertex of one
of the small squares. $P_1$, $P_2$ and $P_3$ are obtained from $P_0$ by 
rotation by an angle of $\frac{\pi}{2}$, $\pi$ and $\frac{3\pi}{2}$, resp.,
around the center of the large square. The two large squares are glued in 
such a way that the canonical map (of degree 2) from the resulting surface
$X_P$ to the torus corresponding to the large square is ramified exactly
over $P_0$, $P_1$, $P_2$ and $P_3$. For the precise description of the
glueing we refer to \cite{WMS}; here we confine ourselves to an example 
with $n=5$:\\[-5mm]
\begin{center}
\begin{picture}(7,7)
\label{bildcase3}
\put(1,1){\framebox(1,1){}}
\put(1,2){\framebox(1,1){}}
\put(1,3){\framebox(1,1){}}
\put(1,4){\framebox(1,1){}}
\put(1,5){\framebox(1,1){}}
\put(2,1){\framebox(1,1){}}
\put(2,2){\framebox(1,1){}}
\put(2,3){\framebox(1,1){}}
\put(2,4){\framebox(1,1){}}
\put(2,5){\framebox(1,1){}}
\put(2,1){\framebox(1,1){}}
\put(2,2){\framebox(1,1){}}
\put(2,3){\framebox(1,1){}}
\put(2,4){\framebox(1,1){}}
\put(2,5){\framebox(1,1){}}
\put(3,1){\framebox(1,1){}}
\put(3,2){\framebox(1,1){}}
\put(3,3){\framebox(1,1){}}
\put(3,4){\framebox(1,1){}}
\put(3,5){\framebox(1,1){}}
\put(4,1){\framebox(1,1){}}
\put(4,2){\framebox(1,1){}}
\put(4,3){\framebox(1,1){}}
\put(4,4){\framebox(1,1){}}
\put(4,5){\framebox(1,1){}}
\put(5,1){\framebox(1,1){}}
\put(5,2){\framebox(1,1){}}
\put(5,3){\framebox(1,1){}}
\put(5,4){\framebox(1,1){}}
\put(5,5){\framebox(1,1){}}

\put(2.05,0.45){$P'_0$}
\put(2.1,6.2){$P'_0$}
\put(6.1,2.1){$P'_1$}
\put(0.2,2.1){$P'_1$}
\put(4.6,6.25){$P'_2$}
\put(4.6,0.4){$P'_2$}
\put(0.3,4.6){$P'_3$}
\put(6.1,4.6){$P'_3$}

\put(2,1){\circle*{.3}}
\put(2,6){\circle*{.3}}
\put(6,2){\circle*{.3}}
\put(1,2){\circle*{.3}}
\put(5,6){\circle*{.3}}
\put(5,1){\circle*{.3}}
\put(1,5){\circle*{.3}}
\put(6,5){\circle*{.3}}

\put(2.1,1){\rule[-0.06\unitlength]{2.8\unitlength}{0.12\unitlength}}
\put(2.1,6){\rule[-0.06\unitlength]{2.8\unitlength}{0.12\unitlength}}
\put(5.96,2.1){\rule{0.12\unitlength}{2.8\unitlength}}
\put(0.96,2.1){\rule{0.12\unitlength}{2.8\unitlength}}

\put(3.4,.65){$a$}
\put(6.2,3.4){$b$}
\put(3.4,6.2){$c$}
\put(0.6,3.4){$d$}

\put(0.2,.6){$O_1$}
\put(1,1){\circle*{.1}}
\put(1,1){\circle{.2}}

\put(3.6,3.6){$M_1$}
\put(3.5,3.5){\circle*{.1}}
\put(3.5,3.5){\circle{.2}}

\end{picture}
\hspace*{5mm}
\begin{picture}(7,7)
\label{bildcase3}
\put(1,1){\framebox(1,1){}}
\put(1,2){\framebox(1,1){}}
\put(1,3){\framebox(1,1){}}
\put(1,4){\framebox(1,1){}}
\put(1,5){\framebox(1,1){}}
\put(2,1){\framebox(1,1){}}
\put(2,2){\framebox(1,1){}}
\put(2,3){\framebox(1,1){}}
\put(2,4){\framebox(1,1){}}
\put(2,5){\framebox(1,1){}}
\put(2,1){\framebox(1,1){}}
\put(2,2){\framebox(1,1){}}
\put(2,3){\framebox(1,1){}}
\put(2,4){\framebox(1,1){}}
\put(2,5){\framebox(1,1){}}
\put(3,1){\framebox(1,1){}}
\put(3,2){\framebox(1,1){}}
\put(3,3){\framebox(1,1){}}
\put(3,4){\framebox(1,1){}}
\put(3,5){\framebox(1,1){}}
\put(4,1){\framebox(1,1){}}
\put(4,2){\framebox(1,1){}}
\put(4,3){\framebox(1,1){}}
\put(4,4){\framebox(1,1){}}
\put(4,5){\framebox(1,1){}}
\put(5,1){\framebox(1,1){}}
\put(5,2){\framebox(1,1){}}
\put(5,3){\framebox(1,1){}}
\put(5,4){\framebox(1,1){}}
\put(5,5){\framebox(1,1){}}
\put(2.05,0.45){$P'_0$}
\put(2.1,6.2){$P'_0$}
\put(6.1,2.1){$P'_1$}
\put(0.2,2.1){$P'_1$}
\put(4.6,6.25){$P'_2$}
\put(4.6,0.4){$P'_2$}
\put(0.3,4.6){$P'_3$}
\put(6.1,4.6){$P'_3$}

\put(2,1){\circle*{.3}}
\put(2,6){\circle*{.3}}
\put(6,2){\circle*{.3}}
\put(1,2){\circle*{.3}}
\put(5,6){\circle*{.3}}
\put(5,1){\circle*{.3}}
\put(1,5){\circle*{.3}}
\put(6,5){\circle*{.3}}

\put(2.1,1){\rule[-0.06\unitlength]{2.8\unitlength}{0.12\unitlength}}
\put(2.1,6){\rule[-0.06\unitlength]{2.8\unitlength}{0.12\unitlength}}
\put(5.96,2.1){\rule{0.12\unitlength}{2.8\unitlength}}
\put(0.96,2.1){\rule{0.12\unitlength}{2.8\unitlength}}

\put(3.4,.65){$c$}
\put(6.2,3.4){$d$}
\put(3.4,6.2){$a$}
\put(0.6,3.4){$b$}

\put(0.2,.6){$O_2$}
\put(1,1){\circle*{.1}}
\put(1,1){\circle{.2}}

\put(3.6,3.6){$M_2$}
\put(3.5,3.5){\circle*{.1}}
\put(3.5,3.5){\circle{.2}}

\end{picture}
\\
{\em 
\begin{tabular}{l}
Highlighted edges: those with same labels are glued.\\ 
Other edges: opposite edges are glued.
\end{tabular}}
\end{center}
Here, $P'_0$, $P'_1$, $P'_2$ and $P'_3$ are the preimages of
the ramification points $P_0$, $P_1$, $P_2$ and $P_3$,
$M_1$ and $M_2$ the preimages of $M$ and $O_1$, $O_2$ the preimages of $O$, 
respectively. { As in Section 3.5, $O$ and $M$ denote the fixed
points of the automorphism $\bar c$ of order $4$ on the torus $E$.}\\

For $n=3$ there are only two different possibilities for the orbit 
 $P_0$, $P_1$, $P_2$, $P_3$:\\[3mm]
\hspace*{-1cm}
\setlength{\unitlength}{.8cm}
\begin{picture}(4.5,4.5)
\label{bildcase3}
\put(1,1){\line(1,0){3}}
\put(1,2){\line(1,0){3}}
\put(1,3){\line(1,0){3}}
\put(1,4){\line(1,0){3}}
\put(1,1){\line(0,1){3}}
\put(2,1){\line(0,1){3}}
\put(3,1){\line(0,1){3}}
\put(4,1){\line(0,1){3}}

\put(1.4,1.4){}
\put(2.4,1.4){}
\put(3.4,1.4){}
\put(1.4,2.4){}
\put(2.4,2.4){}
\put(3.4,2.4){}
\put(1.4,3.4){}
\put(2.4,3.4){}
\put(3.4,3.4){}

\put(2,1){\circle*{.3}}
\put(4,2){\circle*{.3}}
\put(3,4){\circle*{.3}}
\put(1,3){\circle*{.3}}

\put(2,1){\rule[-0.06\unitlength]{1\unitlength}{0.12\unitlength}}
\put(2,4){\rule[-0.06\unitlength]{1\unitlength}{0.12\unitlength}}
\put(1,2){\rule{.12\unitlength}{1\unitlength}}
\put(4,2){\rule{0.12\unitlength}{1\unitlength}}

\put(2.35,0.5){$d$}
\put(4.2,2.3){$a$}
\put(2.35,4.15){$c$}
\put(.65,2.3){$b$}
\end{picture}
\begin{picture}(4.5,4.5)
\label{bildcase3}
\put(1,1){\line(1,0){3}}
\put(1,2){\line(1,0){3}}
\put(1,3){\line(1,0){3}}
\put(1,4){\line(1,0){3}}
\put(1,1){\line(0,1){3}}
\put(2,1){\line(0,1){3}}
\put(3,1){\line(0,1){3}}
\put(4,1){\line(0,1){3}}

\put(2,1){\circle*{.3}}
\put(4,2){\circle*{.3}}
\put(3,4){\circle*{.3}}
\put(1,3){\circle*{.3}}

\put(2,1){\rule[-0.06\unitlength]{1\unitlength}{0.12\unitlength}}
\put(2,4){\rule[-0.06\unitlength]{1\unitlength}{0.12\unitlength}}
\put(1,2){\rule{.12\unitlength}{1\unitlength}}
\put(4,2){\rule{0.12\unitlength}{1\unitlength}}
\put(2.35,0.5){$c$}
\put(4.2,2.3){$b$}
\put(2.35,4.15){$d$}
\put(.65,2.3){$a$}

\put(1.3,1.4){}
\put(2.3,1.4){}
\put(3.3,1.4){}
\put(1.3,2.4){}
\put(2.3,2.4){}
\put(3.3,2.4){}
\put(1.3,3.4){}
\put(2.3,3.4){}
\put(3.3,3.4){}

\end{picture}
\hspace*{9mm}
\begin{picture}(4.5,4.5)
\label{bildcase3}
\put(1,1){\line(1,0){3}}
\put(1,2){\line(1,0){3}}
\put(1,3){\line(1,0){3}}
\put(1,4){\line(1,0){3}}
\put(1,1){\line(0,1){3}}
\put(2,1){\line(0,1){3}}
\put(3,1){\line(0,1){3}}
\put(4,1){\line(0,1){3}}

\put(2,2){\circle*{.35}}
\put(2,3){\circle*{.35}}
\put(3,2){\circle*{.35}}
\put(3,3){\circle*{.35}}

\put(3,2){\rule[-0.06\unitlength]{1\unitlength}{0.12\unitlength}}
\put(2,2){\rule[-0.06\unitlength]{1\unitlength}{0.12\unitlength}}
\put(1,2){\rule[-0.06\unitlength]{1\unitlength}{0.12\unitlength}}

\put(2,2){\rule{.12\unitlength}{1\unitlength}}
\put(3,1){\rule{0.12\unitlength}{1\unitlength}}
\put(3,3){\rule{0.12\unitlength}{1\unitlength}}

\put(3.2,1.2){$b$}
\put(2.7,1.2){$a$}
\put(3.5,1.7){$d$}
\put(3.5,2.1){$c$}
\put(2.5,2.1){$f$}
\put(2.5,1.6){$e$}
\put(1.2,2.1){$h$}
\put(1.2,1.65){$g$}
\put(2.2,2.4){$j$}
\put(1.7,2.4){$i$}
\put(3.2,3.4){$l$}
\put(2.7,3.4){$k$}

\end{picture}
\begin{picture}(4.5,4.5)
\label{bildcase3}
\put(1,1){\line(1,0){3}}
\put(1,2){\line(1,0){3}}
\put(1,3){\line(1,0){3}}
\put(1,4){\line(1,0){3}}
\put(1,1){\line(0,1){3}}
\put(2,1){\line(0,1){3}}
\put(3,1){\line(0,1){3}}
\put(4,1){\line(0,1){3}}

\put(2,2){\circle*{.35}}
\put(2,3){\circle*{.35}}
\put(3,2){\circle*{.35}}
\put(3,3){\circle*{.35}}

\put(3,2){\rule[-0.06\unitlength]{1\unitlength}{0.12\unitlength}}
\put(2,2){\rule[-0.06\unitlength]{1\unitlength}{0.12\unitlength}}
\put(1,2){\rule[-0.06\unitlength]{1\unitlength}{0.12\unitlength}}

\put(2,2){\rule{.12\unitlength}{1\unitlength}}
\put(3,1){\rule{0.12\unitlength}{1\unitlength}}
\put(3,3){\rule{0.12\unitlength}{1\unitlength}}

\put(3.2,1.2){$a$}
\put(2.7,1.2){$b$}
\put(3.5,1.7){$c$}
\put(3.5,2.1){$d$}
\put(2.5,2.1){$e$}
\put(2.5,1.6){$f$}
\put(1.2,2.1){$g$}
\put(1.2,1.6){$h$}
\put(2.2,2.4){$i$}
\put(1.7,2.4){$j$}
\put(3.2,3.4){$k$}
\put(2.7,3.4){$l$}

\end{picture}
\vspace*{-7mm}
\begin{center}
{\em 
\begin{tabular}{l}
Highlighted edges: those with same labels are glued.\\ 
Other edges: opposite edges are glued.
\end{tabular}}
\end{center}

In this case
the two corresponding origamis are related by an affine diffeomorphism
which is given explicitly in \cite[Example 3.4]{WMS};
therefore their origami curves are equal by 
Proposition~\ref{gleich}\,b).\\[2mm]
Let $d_P:X\to E$ denote the covering of degree $2n^2$ that defines the 
origami $D_P$. Recall from Section \ref{ocurves} that the points on $C(D_P)$ 
are given by the Riemann surfaces $X_A$ for $A\in\slzwei(\RR)$. Here $X_I$
corresponds to the intersection point of $C(D_P)$ with $C(W)$ described
in Section \ref{cross}. The affine diffeomorphism $c$ on $X_I$ has derivative
$S=\begin{pmatrix}0&-1\\1&0\end{pmatrix}$ and therefore is not an element of
the automorphism group of $D_P$ (cf.~Section \ref{autos}). Explicitly, $c$
defines an affine diffeomorphism $c_A$ on $X_A$ with derivative 
$A\cdot S\cdot A^{-1}$ which, in general, is not holomorphic.\\[2mm]
On the other hand, the square $\tau = c^2$ of $c$ has derivative $-I$
which is central in $\slzwei(\RR)$; therefore $\tau$ defines an affine
and holomorphic automorphism $\tau_A$ of order 2 and derivative $-I$
on each $X_A$. By inspection of the possible
positions of $P$, see \cite[Section~3.2]{WMS}, one finds that, in 
{ any case}, $\tau$ fixes the inverse images of $O$ and $M$ under $d_P$ (and no
other points of $X$).

\subsection{Coverings with given ramification data}
\label{Htilde}

The origamis $D_P$ introduced in the previous section all are defined
by coverings with the same ramification behaviour. It is classically known
that coverings with a prescribed kind of ramification are classified by an 
algebraic variety, called a {\it Hurwitz space}.\\[2mm]
More precisely Hurwitz spaces classify coverings $p:X\to Y$ of compact 
Riemann surfaces with the following data fixed: the degree of $p$, the genus 
of $Y$, the number of critical points and the ramification orders of their
preimages. It is also possible to specify a certain geometric configuration
of the critical points and/or the monodromy of the covering. Note that
once $Y$ and the critical points of $p$ together with their ramification
orders are fixed, there are only finitely many possibilities for the 
monodromy homomorphism $\mu:\pi_1(Y - \{\mbox{critical points}\})\to S_n$,
where $n$ is the degree of $p$. Thus restricting the monodromy typically
leads to an irreducible component of the Hurwitz space defined by the other 
data.\\[2mm]
In our case, we have coverings of degree 2 of an elliptic curve, that
are ramified over 4 points $P$, $P'$, $Q$, $Q'$. Moreover, for a suitable
choice of a base point on the elliptic curve, we have $P' = -P$ and
$Q' = -Q$. Note that we can always achieve $P' = -P$ by choosing an
appropriate base point; but then ``$Q' = -Q$'' is a proper condition.
We cannot formulate in an algebraic way the condition that
$P_0$, $P_1$, $P_2$, $P_3$ form an orbit under the affine diffeomorphism
$\bar c_A$ induced by $\bar c$,
since we saw above that $\bar c_A$ is in general not an algebraic automorphism 
on the elliptic curve $E_A$.\\[2mm]
We denote by $\tilde H$ the Hurwitz space of isomorphism classes of such 
coverings. Thus every element of $\tilde H$ is represented by a pair 
$(X,p)$ where $X$ is a compact Riemann surface and $p:X\to E$ is a covering
of degree 2 to an elliptic curve $E$ with a base point $O$ such that the
four critical points $P$, $Q$, $P'$, $Q'$ of $p$ satisfy $P'=-P$ and
$Q'=-Q$. 
Two such coverings $p:X\to E$ and
$p':X'\to E'$ are considered equivalent (or isomorphic) if there are
isomorphisms $\varphi:X\to X'$ and $\bar\varphi:E\to E'$ such that 
$p'\circ\varphi = \bar\varphi\circ p$.\\[2mm]
By the Riemann-Hurwitz formula, for every $(X,p)\in\tilde H$, the genus of 
$X$ is 3. Therefore there is a natural morphism $\pi$ (of algebraic varieties)
from $\tilde H$ to the moduli space $M_3$ of isomorphism classes of Riemann 
surfaces of genus 3.
\begin{prop}
\label{V4}
Let $(X,p)$ be a point in $\tilde H$. Then
$\emAut(X)$ contains a subgroup isomorphic to the Klein four group $V_4$.
\end{prop}
\begin{proof}
Every covering of degree 2 is normal, therefore there is an automorphism
$\sigma$ of $X$ that interchanges the two inverse images under $p$; the
ramification points $\tilde P$, $\tilde P'$, $\tilde Q$, $\tilde Q'$, 
that are mapped to $P$, $P'$, $Q$, and $Q'$ by $p$, are the fixed points of
$\sigma$. \\
Next we want to show that the automorphism $[-1]$ on the elliptic curve
$E$ lifts to an automorphism $\tau$ on $X$ of order 2 that commutes with 
$\sigma$. Note that $[-1]$ also is an automorphism of 
$E^{**} := E-\{P,P',Q,Q'\}$
and thus induces an automorphism $\tilde \tau$ of $\pi_1(E^{**})$. It follows
from surface topology that $[-1]$ lifts to $X$ if and only if $\tilde\tau$
preserves the subgroup $U$ of $\pi_1(E^{**})$ that corresponds to the unramified
covering $p:X-\{\tilde P, \tilde P', \tilde Q, \tilde Q'\}\to E^{**}$.\\
We take a fixed point $Z_0$ of $[-1]$ as base point for $\pi_1(E^{**})$.
Let $l_x$ and $l_y$ be simple geodesic loops around $Z_0$ that are mapped 
onto their inverses by $[-1]$ and such that their homotopy classes $x$ and $y$
generate
$\pi_1(E)$. Together with the loops $l_P$, $l_{P'}$, $l_Q$ and $l_{Q'}$
around the critical points, we obtain a set of generators for $\pi_1(E^{**})$.
Using them, the automorphism $\tilde\tau$ of order 2 is given by 
\begin{equation}
\tilde\tau(x) = x^{-1},\quad \tilde\tau(y)=y^{-1},\quad \tilde\tau(l_P)=l_{P'},
\quad\tilde\tau(l_Q)=l_{Q'}
\end{equation}
The subgroup $U$ is the kernel of the monodromy homomorphism 
$\mu=\mu_P:\pi_1(E^{**})\to S_2$. Since $p$ is ramified over $P$, $P'$,
$Q$ and $Q'$, we have $\mu(l_P)=\mu(l_{P'})=\mu(l_Q)=\mu(l_{Q'})=(1\ 2)$.
Moreover $\mu(x)=\mu(x^{-1})$ and $\mu(y)=\mu(y^{-1})$, thus the kernel
of $\mu$ clearly is preserved by $\tilde\tau$. Hence $\tilde\tau$ induces
an automorphism $\tau$ of $X$ of order 2 which by construction acts on the
fibers of $p$ and thus commutes with $\sigma$. 
\end{proof}
Note that a covering $p:X\to E$ of degree 2 is completely determined by 
its deck transformation $\sigma$. Therfore we can describe the Hurwitz
space $\tilde H$ also as the set of equivalence classes of pairs $(X,\sigma)$,
where $X$ is a Riemann surface of genus 3 and $\sigma\in\Aut(X)$ an involution
with exactly 4 fixed points whose images on $E=X/<\sigma>$ are symmetric
w.\,r.\,t.\ some point on $E$. According to the proof of Proposition \ref{V4},
the last condition is equivalent to the existence of a second involution $\tau$
which commutes with $\sigma$.\\[2mm]
The equivalence relation can be rewritten as
\begin{equation}
\label{equiv}
(X,\sigma)\cong(X',\sigma')\Leftrightarrow \exists \mbox{ isomorphism }\varphi:
X\to X'\mbox{ s.\,t. }\sigma'=\varphi\sigma\varphi^{-1}
\end{equation}

\subsection{Origami curves in {\boldmath $H$}}
\label{dense}
If we are given an elliptic curve $E$ and four points
$P$, $-P$, $Q$, $-Q$ on $E$, there are four different coverings of $E$
of degree 2 that are ramified exactly over these points: this is due to the
fact, pointed out at the end of the proof of Proposition \ref{V4}, that the
monodromy homomorphism $\mu:\pi_1(E^{**})\to S_2$ is fixed by these data
except for the choice of $\mu(x)$ and $\mu(y)$.
The covering surface $X$ is hyperelliptic
if and only if it admits an involution with 8 fixed points. A (general) point
$(X,\sigma)\in\tilde H$ has 3 involutions: $\sigma$, $\tau$ and $\sigma\tau$
in the notation of Proposition \ref{V4}. By definition, $\sigma$ has 4 fixed
points. Possible fixed points of $\tau$ and $\sigma\tau$ are the inverse
images of the four fixed points of $[-1]$ on $E$.
If $\tau$ fixes exactly 4 of them, the other four are fixed by $\sigma\tau$; 
if $\tau$ has 8 fixed points, $\sigma\tau$ has none, and vice versa. In 
particular, at least one of $\tau$ and $\sigma\tau$ has fixed points on $X$;
thus we may assume that $\tau$ has a fixed point.\\[2mm]
It is easy to see that for precisely one choice of the monodromy map,
$\tau$ has 8 fixed points. 
For the other three choices of $\mu$, {$\tau$ and $\sigma\tau$ have 
4 fixed points each.
\begin{defn}
Let $H$ be the subspace of $\tilde H$ that consists of the pairs
$(X,\sigma)$ for which the involution $\tau$ has exactly 4 fixed points.
\end{defn}

For $(X,\sigma)\in H$, the quotients of $X$ by the involutions $\sigma$,
$\tau$ and $\sigma\tau$ all have genus 1. We shall see in Corollary \ref{G0}
that for generic $(X,\sigma)\in H$, $\Aut(X)=\{\id, \sigma, \tau, \sigma\tau\}$;
in particular, $\Aut(X)$ does not contain an involution with quotient curve
isomorphic to $\PP^1$, i.\,e.\ $X$ is not hyperelliptic. Thus we could
alternatively characterize $H$ as the closure in $\tilde H$ of the set of
those $(X,\sigma)$ for which $X$ is not hyperelliptic.}\\[3mm]
{\bf Remark.}
For every $n$-torsion point $P$ on $E_{-1}$, the origami $D_P$ defines
an algebraic curve $C_P$ in $H$. { In the same way, $W$ defines 
an algebraic curve $C_W$ in $H$.}\\[2mm]
{ Namely, $\pi:\tilde H\to M_3$ is a finite morphism of varieties, 
and the origami curves $C(D_P)$ resp.\ $C(W)$ are algebraic curves in 
$M_3$ which lie in the image $\pi(\tilde H)$. This was shown in 
Section 4.2 for $D_P$ and in Section 3.5 for $W$. The inverse image 
$C_P:=\pi^{-1}(C(D_P))$ resp.\ $C_W = \pi^{-1}(C(W))$ is therefore an algebraic
curve in $\tilde H$. 
It remains to show that they lie in $H$. But  for $C_P$, this follows from 
the observation in Section \ref{dp} that the automorphism
$\tau$ has 4 fixed points on every Riemann surface $X_A$ corresponding to
a point $(X_A,\sigma)\in C_P$. On $W$, $\tau$ is the automorphism $c^2$
(or $-1\in Q$) described in Section 3.2 which fixes the four vertices
of the squares.}
\begin{thm}
\label{dense}
The union of the curves $C_P$, $P$ a torsion point on $E_{-1}$, is dense 
in $H$.
\end{thm}
\begin{proof}
We shall prove the statement for the complex topology on the Hurwitz space
$H$. Then it holds a fortiori also for the Zariski topology on the algebraic
variety $H$.\\[1mm]
Let $(X,\sigma)$ be a point in $H$ with non-hyperelliptic $X$. This point is 
determined by an elliptic curve $E$ (with a base point $O$), two pairs 
$(P,-P)$ and $(Q,-Q)$ of opposite points on $E$, and a choice
of one of three possible monodromy maps. Clearly we can approximate $P$
and $Q$ by torsion points, more precisely by pairs $P_n$, $Q_n$ of $n$-torsion
points for the same $n$. For a suitable choice of the monodromy map, 
the twofold covering 
of $E$ ramified over $P_n$, $-P_n$, $Q_n$ and $-Q_n$ then defines a point
$(X_n,\sigma_n)\in H$ which is close to $(X,\sigma)$. Note that the 
composition of the covering $X_n\to E$ of degree 2 with multiplication by $n$ 
on $E$ is an origami $O$, hence $(X_n,\sigma_n)$ lies on the inverse image
$C = \pi^{-1}(C(O))$ of the origami curve
$C(O)$ in $M_3$. In the sequel, we also call $C \subset H$ an origami curve.\\[1mm]
We now have to show that there is an $n$-torsion point $P$ on $E_{-1}$
such that $(X_n,\sigma_n)\in C_P$. We shall in fact prove the stronger
statement $C=C_P$, but only for the case that $n=p$ is a prime (which
clearly suffices to approximate the original point $(X,\sigma)$).\\[1mm]
Recall from Proposition \ref{gleich} that the origami curves $C$ and $C_P$ are 
equal if and only if the corresponding subgroups $U$ and $U'$ of $F_2$ 
are mapped onto each other by an automorphism of $F_2$. If  such
an automorphism exists, it is induced by an affine diffeomorphism 
$\bar f:E\to E_{-1}$
which transforms the configuration $\{P_p,Q_p,-P_p,-Q_p\}$ of $p$-torsion 
points on $E$ to an orbit $\{P,\bar c(P),-P,-\bar c(P)\}$ under $\bar c$
on $E_{-1}$ (for some $P$),
and which moreover respects the monodromy. 
We shall see in the proof of Proposition \ref{indexprop} that, 
for a given pair of $p$-torsion points on $E$, all three non-hyperelliptic 
choices of the monodromy map lead to points on the same origami curve; 
so we do not have to care about the monodromy. \\[1mm]
To find such an affine diffeomorphism $\bar f$,
consider a point $(X',\sigma')$ on the origami curve $C$
through $(X_p,\sigma_p)$ which lies over the elliptic curve $E_{-1}$.
To $(X',\sigma')$ there corresponds an affine transformation $A\in\slzwei(\RR)$
that maps $E$ into $E_{-1}$, and $P_p$ and $Q_p$ to $p$-torsion points on 
$E_{-1}$. We identify the group $E_{-1}[p]$ of $p$-torsion points on $E_{-1}$ 
with
$(\ZZ/p\ZZ)^2$ and keep the notation $P_p$ and $Q_p$ for the points in 
$(\ZZ/p\ZZ)^2$ induced by the original $p$-torsion points on $E$.
Thus in order to show that the origami curves $C$ (defined by $O$) and 
$C_P$ (defined by $P\in E_{-1}[p]$) are equal it suffices to find some
$B\in\slzwei(\ZZ/p\ZZ)$ and some $P\in(\ZZ/p\ZZ)^2$ such that
\begin{equation}
\label{matrixB}
\{B(P_p),B(Q_p)\} = \{P,\bar c(P)\}.
\end{equation}
Since $\bar c$ is the rotation by $\frac{\pi}{2}$, it acts on $E_{-1}[p]$
through the matrix $ S=\begin{pmatrix}\scs{0}&\scs-1\\\scs1&\scs0\end{pmatrix}$.
\\[1mm]
Now we identify $\ZZ/p\ZZ$ with the field $\FF_p$ and use linear algebra 
over $\FF_p$: we consider $P_p$ and $Q_p$ as vectors $\begin{pmatrix}\scs a\\
\scs b\end{pmatrix}$, $\begin{pmatrix}\scs c\\\scs d\end{pmatrix}$ in $\FF_p^2$;
we may assume that $\det\begin{pmatrix}\scs a&\scs c\\\scs b&\scs d\end{pmatrix}
\not=0$, because we can approximate $(X,\sigma)$ by $p$-torsion points 
that are not collinear. Choose $M\in\slzwei(\FF_p)$ with $M(P_p)=
\begin{pmatrix}\scs 1\\\scs 0\end{pmatrix}$, and let 
$\begin{pmatrix}\scs u\\\scs v\end{pmatrix} = M(Q_p)$. Since $v=ad-bc$ is
invertible in $\FF_p$, we find a matrix $S''\in\slzwei(\FF_p)$
with first column $\begin{pmatrix}\scs u\\\scs v\end{pmatrix}$ and 
$\mbox{tr}(S'')=0$. For the conjugate $S':=M^{-1}S''M$ we then obtain
$S'(P_p)=Q_p$.\\
We show in Lemma \ref{conj} below that $S'$ is
conjugate to $S$ or to $S^{-1}$ in $\slzwei(\FF_p)$, say $S' = B^{-1}SB$.
Then
\begin{eqnarray*}
SB(P_p)=BS'(P_p)=B(Q_p),
\end{eqnarray*}
i.\,e.\ (\ref{matrixB}) is satisfied for $P=B(P_p)$, and the proposition is proved.
\end{proof}
\begin{lemma}
\label{conj}
Let $p$ be a prime and $T\in\emslzwei(\FF_p)$ with $\mbox{\em tr}(T)=0$.
Then $T$ is conjugate in $\emslzwei(\FF_p)$ to $S$ or to $S^{-1}$, where
$S=\begin{pmatrix}0&-1\\1&0\end{pmatrix}$.
\end{lemma}
\begin{proof}
First observe that $T^2=-I$ and the characteristic polynomial of $T$ is
$x^2+1$. We distinguish two cases:\\[1mm]
1. If $p\equiv-1\mbox{ mod }4$, $T$ has no eigenvalue in $\FF_p$. Let $v_1$
be any nonzero vector in $\FF_p^2$ and $v_2=Tv_1$. With respect to the
basis $v_1$, $v_2$, $T$ is represented by the matrix $S$. Hence $T$ and
$S$ are conjugate by an invertible matrix $M\in\mbox{GL}_2(\FF_p)$. Replacing
$v_1$, $v_2$ by $\lambda v_1$ and $\lambda v_2$ for some 
$\lambda\in\FF_p^\times$ multiplies $\det(M)$ by $\lambda^2$. Thus if the
determinant of the original matrix $M$ was a square in $\FF_p^\times$,
we can find a suitable $\lambda$ such that the base change matrix has
determinant 1. Otherwise $-\det(M)$ is a square (since $p\equiv-1\mbox{ mod }4$
and hence $-1$ is not a square in $\FF_p$). We now replace $v_2$ by $-v_2$
and obtain a base change matrix with determinant $-\det(M)$. $T$ is conjugate
to $S^{-1}$ in this case.\\[1mm]
2. If $p\equiv 1\mbox{ mod }4$, $T$ has two eigenvalues $\alpha$ and
$-\alpha$ in $\FF_p$ (satisfying $\alpha^2=-1$). Thus $T$ is conjugate to
$\tilde S=\begin{pmatrix}\alpha&0\\0&-\alpha\end{pmatrix}$ by a matrix $M$ 
whose rows are eigenvectors of $T$. Multiplying the first row of $M$ by
$\det(M)^{-1}$ then yields a base change matrix of determinant 1. 
In particular, $S$ itself is conjugate to $\tilde S$ in $\slzwei(\FF_p)$,
thus all other matrices $T\in\slzwei(\FF_p)$ with $\mbox{tr}(T)=0$ are
conjugate to $S$.
\end{proof}

\subsection{Affine coordinates for {\boldmath $H$}}
\label{U}
\begin{prop}
\label{gleichung}
For any point $(X,\sigma)\in H$, $X$ can be represented by a plane quartic 
with equation
\begin{equation}
\label{abc}
x^4 + y^4 + z^4 + 2ax^2y^2 + 2bx^2z^2 + 2cy^2z^2 = 0
\end{equation}
for some complex numbers $a$, $b$ and $c$.
\end{prop}
\begin{proof}
Since $X$ is not hyperelliptic, the canonical map on $X$ gives an embedding
into the projective plane as a smooth quartic. Moreover, every automorphism
of $X$ is induced by a projective automorphism of $\PP^2$. Thus we can
choose coordinates on $\PP^2$ such that $\sigma$ acts by 
$(x:y:z)\mapsto (-x:y:z)$ and $\tau$ by $(x:y:z)\mapsto (x:-y:z)$. Then 
$\sigma\tau$ acts by $(x:y:z)\mapsto (-x:-y:z) = (x:y:-z)$. A quartic
that is invariant under $\sigma$ and $\tau$ must therefore be a polynomial
in $x^2$, $y^2$ and $z^2$. We can still multiply $x$, $y$ and $z$ by
suitable constants and thus obtain that the coefficients of $x^4$, $y^4$
and $z^4$ are 1, which gives us an equation of the form (\ref{abc}).
\end{proof}
The equation (\ref{abc}) describes a family of plane projective curves
of degree 4 over the affine parameter space $\AAA^3(\CC) = \CC^3$. The 
nonsingular curves in this family are determined by
\begin{prop}
{\bf\em a)} The plane curve $C_{abc}$ with equation $(\ref{abc})$ is singular
if and only if $a=\pm1$ or $b=\pm1$ or $c=\pm1$ or 
$a^2+b^2+c^2-2abc-1=0$.\\[1mm]
{\bf\em b)} Let $V\subset\AAA^3(\CC)$ be the zero set of the polynomial
$(a^2-1)(b^2-1)(c^2-1)(a^2+b^2+c^2-2abc-1)$, $U=\AAA^3(\CC)-V$ and
\[{\cal C} = \{((x:y:z),(a,b,c))\in\PP^2(\CC)\times U:(x:y:z)\in C_{abc}\}.\]
Then the projection $\mbox{\em proj}_2:{\cal C}\to U$ to the second factor
is a family of nonsingular projective curves of genus $3$.
\end{prop}
\begin{proof}
The first part is easily checked using Jacobi's criterion 
$\frac{\partial f}{\partial x}(P)=\frac{\partial f}{\partial y}(P)=
\frac{\partial f}{\partial z}(P)=0$ for a singular point $P$ on the plane
projective curve with equation $f=0$.
The second statement is a consequence of the first.
\end{proof}
It follows from the property of $M_3$ as a coarse moduli space, that the map
$m:U\to M_3$ that sends a point $(a,b,c)$ to the isomorphism class of
$C_{abc}$, is an algebraic morphism. To better understand this map we 
first determine the automorphisms of $U$ and of ${\cal C}$. 
\begin{prop}
\label{AutC}
{\bf\em a)} For every $(a,b,c)\in U$, $C_{abc}$ admits the automorphisms
\[\alpha:(x:y:z)\mapsto (-x:y:z)\ \mbox{ and }\ \beta:(x:y:z)\mapsto (x:-y:z);\]
they generate a Klein four group $G_0$ in the group 
\[\emAut({\cal C}) = 
\{\varphi:{\cal C}\to{\cal C}:\varphi \mbox{ \em automorphism},  
\mbox{ \em proj}_2\circ\varphi = \mbox{ \em proj}_2\}\] 
of relative automorphisms of the family ${\cal C}$.\\[1mm]
{\bf\em b)} The group $L$ of linear automorphisms of $U$ is generated by
the permutations of $a$, $b$ and $c$, and by the map $s:(a,b,c)\mapsto 
(-a,-b,c)$; $L$ is isomorphic to $S_4$.\\[1mm]
{\bf\em c)} There is a short exact sequence
\begin{equation}
\label{exact}
1 \to G_0 \to \emAut({\cal C}) \to L \to 1.
\end{equation}
{\bf\em d)} $\emAut({\cal C})$ has 96 elements and is isomorphic to the
automorphism group of the Fermat curve $C_{000}: x^4+y^4+z^4=0$.
\end{prop}
\begin{proof}
a) is obvious.\\ 
b) Clearly { the group generated by $s$ and the permutations of
$a$, $b$ and $c$} is contained in $\Aut(U)$, and it is easy to check that 
{ it} is isomorphic to $S_4$.\\
On the other hand, any linear automorphism $f$ of $U$ must permute the six 
planes $a=\pm1$, $b=\pm1$
and $c=\pm1$. The nonlinear component of $V$ also must be preserved, which
gives a restriction on the distribution of signs: the product of the three
signs must be $+1$. Thus $f$ is in the group generated by $S_3$ and $s$.\\
c) Observe that every element of $L$ can be lifted to an automorphism of
${\cal C}$: For a permutation of $a$, $b$ and $c$, apply the same permutation
to $z$, $y$ and $x$ (in this order!); $s$ can be lifted to 
$\tilde s:((x:y:z),(a,b,c))\mapsto((ix:y:z),(-a,-b,c))$. The subgroup $G$
of $\Aut({\cal C})$ generated by $G_0$ and these lifts clearly fits into
the exact sequence (\ref{exact}). We shall see later that $G=\Aut({\cal C})$
(see the remark following Proposition \ref{birat}).\\
d) $G$ has $|G_0|\cdot|L|=96$ elements, and all of them map $C_{000}$ onto
itself. { 
Recall that $C_{000}$ is the 
Fermat curve of degree 4 and has precisely 96 automorphisms as we have 
mentioned in Section 3.3.}
\end{proof}
\begin{cor}
\label{G0}
The kernel of the homomorphism $\rho:\emAut({\cal C}) \to \emAut(U)$ is $G_0$.
In particular, for general $(a,b,c)\in U$, $\emAut(C_{abc}) = G_0$.
\end{cor}
\begin{proof}
As observed in the proof of Proposition { \ref{gleichung}},
any automorphism of a curve in ${\cal C}$ is
induced by an automorphism of $\PP^2(\CC)$. Thus any 
$\varphi\in\mbox{ker}(\rho)$ is a linear change of the homogeneous coordinates
$x$, $y$, $z$ that preserves the terms $ax^2y^2$, $bx^2z^2$ and $cy^2z^2$
for all $(a,b,c)\in U$. This is only possible if $\varphi\in G_0$. 
Alternatively, the result can be checked by comparison with the list of 
automorphism groups in genus 3 in \cite{KK}.
\end{proof}

\subsection{Maps to moduli space}
\label{24}
In this subsection we study the relations between the spaces $U$, $H$ and 
$M_3$. To do so, we shall factor the morphism $m:U\to M_3$, $(a,b,c)\mapsto
[C_{abc}]$, in two ways.\\[2mm]
We saw in the proof of Proposition \ref{AutC} that the subgroup $L$
of $\Aut(U)$ is contained in the image of the natural homomorphism
$\rho:\Aut({\cal C})\to \Aut(U)$. Therefore $m$ factors through $U/L$.\\[2mm]
On the other hand we know that , for every $(a,b,c)\in U$, $C_{abc}$ admits 
the automorphism $\alpha:(x:y:z)\mapsto(-x:y:z)$. We find:
\begin{prop}
\label{grad1}
$(a,b,c)\mapsto(C_{abc},\alpha)$ is a surjective morphism $h:U\to H$.
\end{prop}
\begin{proof}
The fixed points of $\alpha$ in $\PP^2(\CC)$ are the point $(1:0:0)$ (which
is not on $C_{abc}$ for any $(a,b,c)\in U$) and all points having $x=0$.
On $C_{abc}$ there are precisely 4 points $(0:y_i:1)$ of this form, where
the $y_i$ are the solutions of $y^4+2cy^2+1=0$ (note that they are all
different because $c^2\not=1$). Thus $E_{abc}:=C_{abc}/<\!\alpha\!>$ is a
curve of genus 1 by the Riemann-Hurwitz formula. 
To see that $(C_{abc},\alpha)$ defines a point in $H$,
it remains to show that the critical points of the covering 
$p:C_{abc}\to E_{abc}$ are symmetric w.\,r.\,t.\ an involution $[-1]$.
For this consider the automorphism $\beta:(x:y:z)\mapsto(x:-y:z)$ of $C_{abc}$
which, like $\alpha$, is an element of $G_0$: $\beta$ also has 4 fixed points 
on $C_{abc}$, different from those of $\alpha$, and it maps every fixed 
point of $\alpha$ to a fixed point of
$\alpha$. Since $\alpha$ and $\beta$ commute, $\beta$ descends to an 
involution on $E_{abc}$ with 4 fixed points that acts as $[-1]$ on the 
critical points of $p$.\\
This defines the morphism $h$; the surjectivity of $h$ was proved in 
Proposition \ref{gleichung}.
\end{proof}
So far we have found a commutative diagram
\begin{center}
$\xymatrix@=6ex{\ar[r]^h\ar[d]_{\tilde q}\ar[dr]_m U&\ar[d]^{\pi}H\\
\ar[r]_q U/L&M_3}$
\end{center}
The final goal in this section is to show
\begin{prop}
\label{birat}
$q$ is birational.
\end{prop}
{\bf Remark.}
From this it follows in particular that $L=\rho(\Aut({\cal C}))$ (because
$m$ also factors through $U/\rho(\Aut({\cal C}))$. Together with Corollary 
\ref{G0} this shows that $\Aut({\cal C})$ fits in the middle of the exact
sequence (\ref{exact}), and thereby finishes the proof of Proposition 
\ref{AutC}.
\begin{proof}[Proof of Proposition \ref{birat}.]
Since $U$ is irreducible, the same holds for $H$, $U/L$ and their image
$m(U)\subset M_3$. Moreover all morphisms in the diagram are finite, so we 
can compare their degrees.\\
Obviously the degree of $\tilde q$ is $|L|=24$.\\
Next we observe that $\mbox{deg}(\pi)=3$: Any covering $C\to E$ of degree 2 is 
induced by an automorphism of order 2 on $C$. Since for a general point
$(C,p)\in H$, $C$ has precisely 3 such involutions, the degree of $\pi$ is 
at most 3. But if some $(C,\alpha)\in H$ would be isomorphic to $(C,\beta)$
for different automorphisms $\alpha$ and $\beta$, these would by (\ref{equiv}) 
have to be 
conjugate in $\Aut(C)$; this is impossible, again because $\Aut(C)=G_0$
for generic $C$ by Corollary \ref{G0}.\\
Finally we claim that $h$ is of degree 8. To see this we have to determine,
for generic $(a,b,c)\in U$, all $(a',b',c')$ such that there is an
isomorphism $\varphi:C_{abc}\to C_{a'b'c'}$ satisfying 
$\varphi\alpha\varphi^{-1} = \alpha$ (see (\ref{equiv})). 
As before, $\varphi$ has to be induced
by an automorphism of $\PP^2(\CC)$ and, since we want $(a,b,c)$ to be generic,
even by an automorphism of ${\cal C}$.\\ 
As a projective transformation, $\varphi$ commutes with $\alpha$ if and only
if it is given by a matrix of the form
\begin{equation}
D_\varphi = \begin{pmatrix}v&0&0\\0&r&s\\0&t&u\end{pmatrix}\in\mbox{PGL}_3(\CC).
\end{equation}
The condition that $\varphi$ maps every curve $C_{abc}$ to a curve of the same
type implies by a straightforward (but unpleasant) calculation that 
\begin{eqnarray*}
v^4=1\quad\mbox{and either}\quad &t=s=0,& r^4=u^4=1\\
\qquad\mbox{or}\quad&
u=r=0,& t^4=s^4=1
\end{eqnarray*} 
This give $4\cdot4\cdot4\cdot2$ matrices; multiplying each entry
by a power of $i$ gives the same element in $\mbox{PGL}_3(\CC)$.
Thus we have found a subgroup of order 32 of $\Aut({\cal C})$ that
generically determines $h$. It contains the kernel $G_0$, therefore its
image $L_H$ in $L$ is a subgroup of order 8.
\end{proof}
{\bf Note:} $q$ is not an isomorphism, nor is $H$ isomorphic to $U/L_H$.\\
This can be seen e.\,g.\ by looking at the Fermat curve $C_{000}$: it is
mapped isomorphically onto $C_{0,3,0}$ by the transformation
$(x:y:z)\mapsto (x+z:\sqrt[4]{8}\,y:x-z)$, which does not extend to an 
automorphism
of ${\cal C}$ since we have seen in Proposition \ref{AutC}\,d) that they 
all fix $C_{000}$.


{\color{blue} }

%% file: vg.tex
\section{The Veech groups to the origamis $D_P$}
\label{VG}

In this section we finally study the Veech groups of the origamis
$D_P$  
introduced in  Section \ref{dp} (with $P$ an $n$-torsion point of $\el$).  
{ Recall that a subgroup of $\slzwei(\ZZ)$ is called 
{\em congruence group of level $n$}, if it contains the 
principal congruence group $\Gamma(n) = \{A \in \slzwei(\ZZ)|\, A \equiv I
\mod n\}$ and $n$ is minimal with this property. As before $I$ denotes the
identity matrix.}
We show that the Veech group
$\Gamma(D_P)$ is a subgroup of index $3$ of a 
congruence subgroup of level $n$ 
(see Proposition \ref{indexdpcor}). For $n$ odd and $P$ in a somehow general 
position, we write the Veech group down explicitly (see Theorem \ref{vgcalc}).
It is a congruence group of level $2n$.\footnote{ Using similar methods
one can calculate the Veech group for $n$ even and $P$ in general position and 
obtains a congruence group of 
level 
$n$. The proof will be carried out elsewhere.}\\

Recall from Section \ref{cross} and Section \ref{dp}
that the origami $D_P$ is given by a degree $2n^2$ covering 
$p:= d_P:X \to \el$, where $\el$ is the torus endowed with the 
translation structure of $\CC/\Lambda_I$ with $\Lambda_I = \ZZ \oplus \ZZ i$. 
The covering $p$ 
splits into a covering $\kappa: X \to \el$ of degree $2$
and the multiplication by $n$: $[n]:\el \to \el$. $\kappa$ is
the quotient by the translation $\sigma$; it is ramified over four
points $P = P_0$, $P_1$, $P_2$ and $P_3$ on $\el$. They
are symmetric with respect to a rotation $\bar{c}$ of angle $\frac{\pi}{2}$. 
Furthermore, $\bar{c}$ can be lifted to an automorphism $c$ on $X$
that has four fixed points. 
\begin{center}
\input{bild7q5.tex}\label{page}\\
{\em $\bar{c}$ acts on $\el$ as rotation of angle $\pi/2$ with fixed points
     $O$ and $M$,\\ 
     it preserves $\{P_0,P_1,P_2,P_3\}$ and lifts to $c$ on $X$.}
\end{center}
We call $\mu$ the translation structure on 
$X\backslash\{\mbox{ramification points
of }p\}$ that is obtained by lifting the one on $\el = \CC/(\ZZ+\ZZ i)$ via 
$p$.\\

We start from the following two observations (proved in Lemma \ref{fixato}):
\begin{itemize}
\item 
Each diffeomorphism $f:X \to X$ which is affine with respect
to $\mu$ descends via $\kappa$ to an affine diffeomorphism
$\bar{f}: \el \to \el$ which one may write as
\begin{equation}\label{fquer}
\bar{f}: \el \to \el, \;\; \bpm x \\y\epm \mapsto A\bpm x\\y\epm + e
\quad \mbox{ with } A \in \slzwei(\ZZ) \mbox{ and } e \in \ZZ^2.
\end{equation}
In order to make (\ref{fquer}) well defined we have to choose some
point to be $O = (0,0)$. We choose one of the two fixed 
points of $\bar{c}$. Furthermore,  
we identify the $n$-torsion points  with $(\ZZ/n\ZZ)^2$ in the
natural way . Since the four ramification points
are symmetric with respect to $\bar{c}$,
we may write them as $P_0 = (p,q)$, $P_1 = (-q,p)$,
$P_2 = (-p,-q)$ and $P_3 = (q,-p)$ with $p,q \in \ZZ/n\ZZ$.
\item
Since $\bar{f}$ lifts to $f$ via $\kappa$, it has to respect
the ramification points, i.e. 
\begin{equation}\label{that}
\bar{f}(S) = S \quad \mbox{ with } 
               S := \{(p,q),(-q,p), (-p,-q), (q,-p)\} 
\end{equation}
\end{itemize}
We then show in Corollary \ref{new},
that the group of affine diffeomorphisms on $X$ projects by $f \mapsto 
\bar{f}$ to a subgroup of index 3 of the group of  diffeomorphisms on $\el$ 
that  fulfill (\ref{that}).
If $n$ is odd, the affine diffeomorphisms have to fix $O$ (Lemma \ref{fixato}) 
and fulfill an additional condition (Lemma \ref{two}). Using this we
write down the Veech groups explicitly for $n$ odd and $S$ in
``general position''.\\

We will use the characterization of origami Veech groups
in Theorem B for the proofs of our claims. Therefore
we describe our setting in terms of subgroups of $F_2$.
We consider the 
fundamental groups of the following punctured surfaces:
\[\trstern := \el - \{\infty\},\quad   \trnstern := \el - [n]^{-1}(\infty) 
\; \mbox{ and }\;  X^* := X - p^{-1}(\infty)\]
Observe that $[n]^{-1}(\infty)$ are precisely 
the $n$-torsion points of $\el$.\\
From the morphisms $\kappa,[n]$ and $p = [n]\circ \kappa$,
we get the following embeddings:
\[U := \pi_1(\Xstern) \;\; \subseteq \;\; 
  H_n := \pi_1(\trnstern) \;\; \subseteq \; \;
 \pi_1(\trstern) = F_2\]
As in Theorem B the group $F_2 = F_2(x,y)$ is the free group on the two 
generators $x$ and $y$. Observe that 
\[
H_n = \mbox{kernel}(\proj_n) \quad \mbox{ with } \;\; 
\proj_n: F_2 \to (\ZZ/n\ZZ)^2, \quad
x \mapsto \bpm1\\0\epm, \;\; 
y \mapsto \bpm0\\1\epm
\]
In particular, $H_n$ is a characteristic subgroup of $F_2$.\\

{\bf Description of the \boldmath$n$-torsion points on $\el$ in terms of these
groups:}\\
The $n$-torsion points are precisely the cusps of $\trnstern= \HH/H_n$.
Hence, they correspond to the equivalence classes 
of the elements $w[x,y]w^{-1} \in H_n$, with $w \in F_2$ and
$w_1[x,y]w_1^{-1} \sim w_2[x,y]w_2^{-1}  \; :\ifff \; 
  w_1w_2^{-1} \in H_n \; \ifff \; 
  \proj_n(w_1) = \proj_n(w_2)$.\\
We may choose the base point $B$ of the fundamental group $\pi_1(\HH/H_n)$ 
such that $[x,y]$ is the loop around the cusp $(0,0)$ (see the Figure 
of $\el$ on page \pageref{page}). Then $l_{a,b} := x^ay^b[x,y]y^{-b}x^{-a}$ 
with $a,b \in \{0,\ldots,n-1\}$ is 
the loop around the cusp  $(\bar{a},\bar{b})$, and
these loops form a system of representatives for the
equivalence classes.\\

{\bf Description of \boldmath$U$ by monodromy of \boldmath$\kappa$:}\\
Since $\kappa:X \to \el$ is a morphism of degree $2$, $U$ is a normal
subgroup of $H_n$ of index 2. Thus it is the kernel of a monodromy map 
$\mu: H_n \to S_2$.\\
From the ramification data of $\kappa$ and the correspondence between
the cusps of $\trnstern$ and the equivalence classes of 
loops around cusps it follows that 
\begin{equation}\label{rd}
\mu(l_{a,b}) =  \bpm 1 & 2 \epm 
  \; \iff \;
  (\bar{a},\bar{b}) 
  \; \in \;  S = \{(p,q), (-q,p), (-p,-q), (q,-p)\}.
\end{equation}
Here, $\bar{a}, \bar{b}$ are the images of $a$, $b$ in $\ZZ/n\ZZ$.
In the following we will also write $l_{a,b} \in S$ meaning
that $(\bar{a},\bar{b}) \in S$ and more generally $l \in S$, if $l$ is
a loop equivalent to some $l_{a,b}\in S$, i.e. if $l$ is a loop around
a ramification point.\\
Since $\trnstern$ is a torus with $n^2$ punctures, its fundamental
group $H_n$ is generated by the simple loops $x^n$, $y^n$ and the
$n^2$ loops $l_{a,b}$ around the punctures. (More precisely the 
fundamental group is the free group in $n^2+1$ generators and
we could omit one of the loops around a puncture. But since there
is none better than the other, we take all of them.)\\
There are precisely four maps 
$\mu_1$, $\mu_2$, $\mu_3$, $\mu_4 : H_n \to S_2$ fulfilling (\ref{rd}):
The images of 
the loops around the cusps are given by (\ref{rd}). For the
other two generators $x$ and $y$ one has the following four possibilities
\[
\begin{array}{ll}
\mu_1:\; x^n \mapsto \bpm 1&2 \epm, \;\; y^n \mapsto \bpm 1&2\epm,\qquad
&\mu_2:\; x^n \mapsto  \id , \;\; y^n \mapsto \id, \qquad\\[1mm]
\mu_3:\; x^n \mapsto \bpm 1&2 \epm, \;\; y^n \mapsto \id,\qquad
&\mu_4:\; x^n \mapsto \id, \;\; y^n \mapsto \bpm 1&2\epm
\end{array}\]
They correspond precisely to the four possible coverings ramified over the
four given points in $S$ that we described in Section \ref{Htilde}.
The coverings are the maps
\begin{equation}\label{ui}
\HH/U_i \to \HH/H_n \;\; \mbox{ with } \;  U_i := \mbox{kernel}(\mu_i),
\end{equation}
induced by the inclusions
$U_i \subseteq H_n$.\\

{\bf Affine diffeomorphisms:}\\
Let us choose a universal
covering $u: \HH \to \Xstern$. We obtain the sequence:
\[\HH \; \stackrel{u}{\to}  \; \Xstern  \; \stackrel{\kappa}{\to}  \;  
\el  \; \stackrel{[n]}{\to}  \; \el .\]
Recall from the remark to Theorem B that 
for the group $\aff(\HH,u)$ of
diffeomorphisms that are affine with respect to the translation structure 
lifted from $\Xstern$ via $u$ (which is the same as the one lifted from 
$E^*$ via $[n]\circ\kappa\circ u$),
one has the following correspondence:
\[\aff(\HH,u) \; \cong  \; \autplus(F_2) \quad  \mbox{ by }  \quad 
   *:\tilde{f} \mapsto \tilde{f}^* := (x \mapsto \tilde{f} x\tilde{f}^{-1}) 
                                   \in \autplus(F_2)\] 
An affine diffeomorphism $\tilde{f}$ descends to a surface $\HH/H$ 
($H\subset F_2$)
if and only if $\tilde{f}^*(H) = H$. The Veech group of $\HH/H$ is 
by Theorem B:
\[\Gamma(\HH/H) = \betadach(\stab(H)) \mbox{ with } 
 \stab(H) = \{ \gamma \in \autplus(F_2)|\, \gamma(H) = H\}.\]

Suppose now, that $f$ is an affine diffeomorphism of $X^*$. 
\begin{lemma}\label{fixato}
$f$ descends via $\kappa$ to an affine diffeomorphism 
$\bar{f}$ on $\el$ fixing the set $S$ of ramification points of $\kappa$.
If $n > 1$ is odd and $P = (p,q)$ is a primitive 
$n$-torsion point (i.e. $n$ is the smallest number such that 
$n\cdot P = 0$), then   $\bar{f}$ has a fixed point at the 
$n$-torsion point $O = (0,0)$.
\end{lemma}

\begin{proof}
Let $\tilde{f}$ be a lift of $f$ to $\HH$ via $u$.
Since $H_n$ is characteristic, one has: $\tilde{f}^*(H_n) = H_n$.
Therefore, by the paragraph before Lemma \ref{fixato}, 
it follows that $f$ descends
via $\kappa$ to $\trnstern = \HH/H_n$.
We may extend this map to a diffeomorphism $\bar{f}$ of the whole surface 
$\el$ and write it as
\[\bar{f}: z \mapsto Az + e \;\; 
  \mbox{ with } A \in \slzwei(\ZZ) \mbox{ and } e 
   \in \ZZ^2.\]
$\bar{f}$ acts on the set of cusps, i.e. on the $n$-torsion 
points of $\el$. Since it can be lifted to $X$ via $\kappa$, it
maps the ramification points of $\kappa$ to themselves, 
i.\,e.\ $\bar{f}(S) = S$.\\ 
Let $n$ be odd. We have to show that $e \equiv (0,0)^t \mod n$. We have:
\begin{equation}\label{sumspq}
\bar{f}(\bpm p\\q \epm) + \bar{f}(\bpm -p\\-q \epm) \equiv 
  \bar{f}(\bpm -q\\p \epm) + \bar{f}(\bpm q\\-p \epm) \equiv 2e \quad 
\mbox{mod n}.
\end{equation}
There are three possibilities to subdivide $S$ into two pairs. 
By (\ref{sumspq}) the two sums 
should be equal mod $n$:
\[\bpm p-q\\q+p\epm = \bpm -p+q\\-q-p \epm, \quad
  \bpm p+q\\q-p\epm = \bpm -p-q\\ p-q \epm, \quad 
  \bpm 0\\0 \epm = \bpm 0\\0\epm \;\; \mbox{ in } \ZZ/n\ZZ\]
In the first two cases one obtains: $2(p-q) \equiv 0$ and $2(p+q) \equiv 0$.
Hence, $4p \equiv 0$ and $4q \equiv 0$. But $P = (p,q)$ is a primitive
$n$-torsion point and $n$ is odd, therefore only the third case is possible
and $2e \equiv (0,0)^t$ mod $n$. Since $n$ is odd, $e \equiv (0,0)^t  \mod n$.
\end{proof}

{It follows from Lemma \ref{fixato}} that in order to obtain the
Veech group of $D_P$, we may restrict to affine diffeomorphisms
of $\el$ that fix the set $S$. By the
correspondence between loops around the cusps of $\HH/H_n$ and 
the $n$-torsion points of $\el$, we have:
\begin{equation}\label{fpl}
\begin{array}{l}
\tilde{f} \in \aff(\HH,u) \mbox{ descends via $\kappa\circ u$ to 
an affine diffeomorphism $\bar{f}$ on $E_{-1}$ that fixes $S$}\\
\hspace*{5mm}\; \iff \; \tilde{f}^* \in \tilde{G} :=  \{\gamma \in \autplus(F_2)|\; 
   \gamma(l_{a,b}) \in S  \; \Leftrightarrow \; l_{a,b}  \in S\} 
\end{array}
\end{equation}

\begin{prop}\label{indexprop}
$\emstab(U)$ is a subgroup of $\tilde{G}$ of index $3$.
\end{prop}

\begin{proof}
We first show that $\stab(U) \subseteq \tilde{G}$:\\
Let $\gamma$ be in $\stab(U)$. With (\ref{rd}), it follows for all 
loops around cusps $l_{a,b}$:
\begin{eqnarray*}
l_{a,b} \in S \, &\iff& \, \mu(l_{a,b}) \neq \id \; \iff \; l_{a,b} \notin U  \; 
\stackrel{\gamma(U) = U}{\iff} \;
  \gamma(l_{a,b}) \notin U \; \\
&\iff& \; \mu(\gamma(l_{a,b})) \neq \id 
\;\iff\; \gamma(l_{a,b}) \in S
\end{eqnarray*}
Thus $\gamma \in \tilde{G}$ and $\stab(U) \subseteq \tilde{G}$.
In fact, this was a group theoretical confirmation for the fact that
$\bar{f}$ has to respect the ramification points,
if it lifts via $\kappa$ to $X$.\\[2mm]
In order to obtain the index, we study the action of $\tilde{G}$ 
on the set $\{U_1,U_2,U_3,U_4\}$ given by 
$U_i \mapsto \gamma(U_i)$. By the definition of $\tilde{G}$, it acts on
these four groups, since they are precisely those groups
that have the same monodromy on the loops around cusps as $U$.\\
$\stab(U) \subseteq \tilde{G}$ is the stabilizing group of 
$U \in \{U_1,U_2,U_3,U_4\}$. Therefore the index $[\tilde{G}:\stab(U)]$ 
of $\stab(U)$ in $\tilde{G}$ is equal to the length of the orbit 
$\tilde{G}\cdot U$.
We will see that there are two orbits: one consisting of a single subgroup,
the other one consisting of the remaining three subgroups including $U$.
This then proves the claim.\\
Let us consider the Riemann surfaces 
$\HH/U_1$, $\HH/U_2$, $\HH/U_3$ and $\HH/U_4$. Recall
that they are the four possible degree 2 coverings of $\el = \HH/H_n$
with ramification locus equal to $S$. Recall furthermore from  \ref{Htilde}
that precisely one of them, let us say $\HH/U_j$, is hyperelliptic. 
But for each $\gamma$ the surface
$\HH/\gamma(U_j)$ is again hyperelliptic: By the remark to
Theorem B, $\gamma$ defines
an affine diffeomorphism $\tilde{f}$ on $\HH$ with $\tilde{f}^* = \gamma$.
Then  $\tilde{f}$ descends to an affine diffeomorphism $f: \HH/U_j \to \HH/\gamma(U_j)$.
Let $A \in \slzwei(\ZZ)$ be its derivative. 
We may conjugate the hyperelliptic involution $\tau$ on $\HH/U_j$ with $f$ 
and obtain the affine map $f\tau f^{-1}$ of derivative 
$A\cdot(-I)\cdot A^{-1}= -I$ 
on $\HH/\gamma(U_j)$, 
where $I$ is the identity matrix. 
Since the derivative is $-I$, it is an automorphism. Furthermore, 
the degree and number of fixed points are the same 
as those of $\tau$ on $\HH/U_j$, therefore it is also a hyperelliptic 
involution.\\
But only for $i = j$, $\HH/U_i$ is hyperelliptic. Therefore $\gamma(U_j) = U_j$.\\
It remains to show that {the other three groups} form only one orbit. 
Let $\gamma_1$ and $\gamma_2$ be the automorphisms
\[\gamma_1: F_2 \to F_2, \;\; x \mapsto  xyx^{-1}, \; y \mapsto x^{-1} 
  \quad \mbox{ and }  \quad
  \gamma_2: F_2 \to F_2, \;\; x \mapsto x, \; y \mapsto yx^{n}\]
One can check that they are both in $\tilde{G}$.
A more geometrical argument is given as follows:
Recall that we had chosen the identification of the $n$-torsion points with 
$(\ZZ/n\ZZ)^2$ such that $(0,0)$ is a fixed point of $\bar{c}$ and therefore
also of $[-1] = \bar{c}^2$.\\
By the remark to Theorem B and the fact that $H_n$ is characteristic, 
$\gamma_1$ and $\gamma_2$ define 
affine diffeomorphisms $\bar{f}_1$, resp. $\bar{f}_2$, on 
$\HH/H_n = \trnstern$ with derivative
\[A_1 = \bpm 0&-1\\1&0\epm, \quad \mbox{respectively } \;\; A_2 = \bpm1&n\\0&1\epm,\]
which we may again extend to diffeomorphisms of $\el$.\\
Observe that both $\gamma_1$ and $\gamma_2$ map $[x,y]$ 
to itself, thus $\bar{f}_1$  and $\bar{f}_2$ have
a fixed point in $(0,0)$. It follows that the action of 
$\gamma_1$ and $\gamma_2$
on the 
equivalence classes of loops $w[x,y]w^{-1}$ is equal to the 
action of the derivatives $A_1$ (resp. $A_2$) on $(\ZZ/n\ZZ)^2$.\\
$A_1$ and $A_2$ both  fix $S$ and thus $\gamma_1$ and $\gamma_2$ are in 
$\tilde{G}$ (see (\ref{fpl})).\\
In the following we determine the orbits of $\gamma_1$ and $\gamma_2$
on $M := \{U_1, U_2, U_3, U_4\}$.
We have:
\[\gamma_1(x^n) = xy^nx^{-1}, \;\; \gamma_1(y^n) = x^{-n}
   \quad \mbox{ and }  \quad
  \gamma_2(x^n) = x^n, \;\; \gamma_2(y^n) = (yx^n)^{n}\]
For $i \in \{1,\ldots,4\}$, one obtains: $\mu_i(\gamma_1(y^n)) = \mu_i(x^{-n}) 
\stackrel{\mbox{\fs in} S_2}{=} \mu_i(x^n)$.\\
For the calculation of $\mu_i(\gamma_1(x^n))$, one has to decompose 
$xy^{n}x^{-1}$ into the chosen generators:
$\mu_i(\gamma_1(x^n)) = \mu_i(xy^{n}x^{-1}) = 
\mu_i(l_{0,0}l_{0,1}\ldots l_{0,n-1}y^n)$.\\
Recall that 
$\mu_i(l_{a,b}) = \bpm1&2\epm$ if and only if $(a,b)$ is a ramification point.
But since the four ramification points are symmetric with respect to
$(0,0)$, either zero or two of them are of the form $(0,b)$. Therefore the
monodromy of the loops around cusps in $\mu_i(\gamma_1(x^n))$ adds up to $\id$
and we have: $\mu_i(\gamma_1(x^n)) = \mu_i(y^n)$.\\
From the calculations of $\mu_i(\gamma_1(x^n))$ and $\mu_i(\gamma_1(y^n))$ 
it follows
that $\gamma_1$ acts  in the following way:
\[\gamma_1 :\;\;  U_1 \mapsto U_1,\quad U_2, \mapsto U_2,\quad
             U_3 \mapsto U_4,\quad U_4 \mapsto U_3\]

Let us now study $\gamma_2$:
Obviously one has
$\mu_i(\gamma_2(x^n)) = \mu_i(x^n)$.
For calculating $\mu_i(\gamma_2(y^n))$, one observes similarly as
above  that $(yx^n)^n$ can be written as product, where the factors
are $x^n$, $n$ times  $y^n$, and loops around cusps. Here the sum of the 
monodromies of the cusps depends on the position of $P$, therefore 
one obtains two cases. If the sum is $\id$ one has
\[\mu_i(\gamma_2)(y^n) =  \mu_i(x^n) + n\cdot \mu_i(y^{n})
\stackrel{n \mbox{ \scs odd}}{=}  \mu_i(x^n) + \mu_i(y^n)\]
and the action of $\gamma_2$ on $M$ is given by:
\[\gamma_2 :\;\;  U_1 \mapsto U_3,\quad U_2, \mapsto U_2,\quad
             U_3 \mapsto U_1,\quad U_4 \mapsto U_4\]
If the sum of the monodromies of the loops around cusps in $(yx^n)^n$
is $\bpm 1&2\epm$, then 
\[\mu_i(\gamma_2)(y^n) =  \mu_i(x^n)+ n\cdot \mu_i(y^{n}) + \bpm 1&2 \epm
\stackrel{n \mbox{ \scs odd}}{=}  \mu_i(x^n) + \mu_i(y^n) + \bpm 1&2 \epm\] 
and the action of $\gamma_2$ is given by:
\[\gamma_2 :\;\;  U_1 \mapsto U_1,\quad U_2, \mapsto U_4,\quad
             U_3 \mapsto U_3,\quad U_4 \mapsto U_2.\]
In both cases the action of $\tilde{G}$ has an orbit of 3 elements.
This finishes the proof.
\end{proof}

From the proposition we obtain the following two conclusions.

\begin{cor}\label{new}
The group of affine diffeomorphisms on $\el$ that lift via $\kappa$
to $X$ is a subgroup of index $3$ of the  group of affine diffeomorphisms 
on $\el$ that fix $S$. 
\end{cor}

\begin{proof}
The isomorphism 
$\star: \aff(\HH,u)\stackrel{\sim}{\rightarrow} \autplus(F_2) $ 
induces the following isomorphisms:
\begin{eqnarray}\label{saturday}
\tilde{G} &\cong& \{\tilde{f}\in \aff(\HH,u)|\; \bar{f} 
\mbox{ fixes } S\} \mbox{\ \ and }\nonumber\\
\stab(U) &\cong& \{\tilde{f}\in \aff(\HH,u)|\; \bar{f} \mbox{ lifts to } X\}.
\end{eqnarray}
We denote as before by $\bar{f}$ the affine diffeomorphism on $\el$ to which 
$\tilde{f}$ descends via $\kappa\circ u$.
Observe that, furthermore, the group of deck transformations of  
$\kappa \circ u: \HH \to \el$ corresponds by $\star$ to the group $\mbox{C}_{H_n} 
\subseteq \autplus(F_2)$ 
of conjugations with elements in $H_n$. 
Taking  in 
(\ref{saturday}) the quotient by $\mbox{C}_{H_n}$, 
the claim follows from the proposition.
\end{proof}

\begin{cor}\label{indexdpcor}
$\Gamma(D_P)$ is a subgroup of index $3$ in $\betadach(\tilde{G})$, 
where $\betadach: \emautplus(F_2) \to \emslzwei(\ZZ)$ is as in Theorem B.
If $n$ is odd, then $\betadach(\tilde{G}) = \emstab(S)$, with
\[\emstab(S) = \{A \in \emslzwei(\ZZ)|\, A\cdot S = S\}.\]
\end{cor}

\begin{proof}
The first claim follows from Proposition \ref{indexprop},
Theorem B, and the following commutative diagram:
\[ \xymatrix{
   1 \ar[r]& \mbox{Inn}(F_2) \cong F_2 \ar[r]& \autplus(F_2)  \ar[r]^{\betadach} & 
     \slzwei(\ZZ)   \ar[r] & 1\\
   1 \ar[r]& H_n \ar[r] \ar@{^{(}->}[u] & \tilde{G}  \ar[r] \ar@{^{(}->}[u] &
     \betadach(\tilde{G})  \ar[r] \ar@{^{(}->}[u]& 1\\
   1 \ar[r]& H_n \ar[r] \ar@{=}[u] & \stab(U)  \ar[r] \ar@{^{(}->}[u] &
     \Gamma(D_P)  \ar[r] \ar@{^{(}->}[u]& 1\\
   }
\]
$\tilde{G}$ projects to the group of affine diffeomorphisms 
$\bar{f}: z \mapsto Az+e$ on $\el$ that map $S$ to itself. If $n$ is odd,
 $e = 0$ by Lemma \ref{fixato}. Therefore $A$ fixes $S$ itself.
\end{proof}

In order to obtain a further condition for an
affine diffeomorphism $f$ on $X$, we study how
the affine diffeomorphism $\bar{f}:\el \to \el$ to which
$f$ descends acts  on the
$2$-torsion points of $\el$. The four $2$-torsion points are
$O = (0,0)$, $M = (\frac{n}{2},\frac{n}{2})$, $(\frac{n}{2},0)$ and
$(0,\frac{n}{2})$. If $n$ is odd, then  $O$ is a fixed
point of $\bar{f}$ by Lemma \ref{fixato}. For this condition, we just needed 
that $\bar{f}$ stabilizes $S$.
In the following Lemma we actually have to use that $\bar{f}$ lifts
to $X$.

\begin{lemma}\label{two}
Let $f \in \emaff(X,\mu)$. 
If  $n$ is odd, then  $\bar{f}$ fixes $M$ and we have:
\[\Gamma(D_P) \;\; \subseteq \;\;
    \Gamma_{u,u} \; := \; \{A \in \emslzwei(\ZZ)\:|\;
                  A \equiv \bpm1&0\\0&1 \epm \; \mbox{\em or } \; 
                  A \equiv \bpm0&-1\\1&0\epm \mod 2\}\]
\end{lemma}

\begin{proof}
Recall that $-I = \bar{c}^2$ is the rotation on $\el$ of angle $\pi$
and has fixed points $M$ and $O$. It can be lifted to the automorphism $\tau$
on $X$  for which the two preimages $O_1$ and $O_2$ of $O$ 
and the two preimages $M_1$ and $M_2$ of $M$ are  fixed points 
(see Section \ref{dp}). 
Geometrically $\tau$ is
a rotation of degree $\pi$ of the two squares that form $X$.\\
Let $A$ be the derivative of $f$. Then $\tau^{-1}f^{-1}\tau f$
is an affine diffeomorphism of $X$ with derivative 
$-IA^{-1}(-I)A = I$. Hence, it is in fact a translation.\\
Furthermore, $\bar{f}(O)=O$ by Lemma \ref{fixato}. Therefore
$f$ either fixes $O_1$ and $O_2$ or it interchanges them.
Since they are both fixed points of $\tau$, we get in both cases:
$f^{-1}\tau^{-1}f\tau(O_1) = O_1$.\\
Hence $f^{-1}\tau^{-1}f\tau$ is a translation and has a fixed point
that is not a ramification point; thus: $f^{-1}\tau^{-1}f\tau =\id$.
It follows that $\tau(f(M_1) = f(\tau(M_1)) = f(M_1)$. Hence, $f(M_1)$
is also a fixed point of $\tau$ and thus one of the points
$O_1$, $O_2$, $M_1$ and $M_2$. But $\bar{f}(O) = O$, therefore
$f(M_1) \in \{M_1, M_2\}$ and $\bar{f}(M) = M$. 
This proves the first part of the claim.\\
Let us now write $\bar{f}$ as 
$ \bar{f}: z \mapsto Az + e \stackrel{\fs n \mbox{ \scs odd}}{=} Az$
with $A \in \slzwei(\ZZ)$.
The coordinates of $M$ are $(\frac{n}{2},\frac{n}{2})$.
It follows that
\[\bpm \frac{n}{2}\\[1mm] \frac{n}{2} \epm \;\;\equiv\;\;
  A \cdot \bpm \frac{n}{2}\\[1mm] \frac{n}{2} \epm  \;\;= \;\;
   \frac{n}{2}\cdot \bpm a+b\\[1mm] c+d \epm \mbox{ mod } n
  \quad \mbox{ with } A = \bpm a&b\\c&d\epm
\]
Thus $a+b$ and $c+d$ have to be odd. This is equivalent with $A$ being in $\Gamma_{u,u}$.
\end{proof}

We may now use these conditions 
in order to calculate the Veech groups explicitly if $n$ is odd and 
$P = P_0 = (p,q)$ and $P_1 = (-q,p)$ are in general position, i.e.:\\
\begin{equation}\label{B} 
\bar{B}  = \bpm p&-q\\q&p\epm \mbox{ is invertible over $\ZZ/n\ZZ$}. 
\end{equation} 

\begin{thm} \label{vgcalc}
Let $n$ be odd and $P = (p,q)$ be an $n$-torsion point such
that $P$ and $P_1 = \bar{c}(P) = (q,-p)$ 
are in general position in the sense explained above. Then
\begin{eqnarray*}
\Gamma(D_P) =&&   
   \{A \in \emslzwei(\ZZ)|\, A \equiv \pm \bpm 1&0\\0&1\epm \; \mbox{ \em or } \; 
                           A \equiv \pm \bpm0&-1\\1&0\epm \; \mbox{\em mod } n\}   \\[2mm] 
   &&\cap \quad  \{A \in \emslzwei(\ZZ)|\, A \equiv \pm \bpm 1&0\\0&1\epm \; 
                                              \mbox{ \em or } \; 
                           A \equiv \pm \bpm0&-1\\1&0\epm \; \mbox{\em mod } 2\}\\
           =&&
   \begin{array}[t]{l}
   \{A = \bpm a &b \\c&d\epm \in \emslzwei(\ZZ)|\,
           a+c \mbox{ \em and } b+d \mbox{ \em odd } \mbox{\em and } \\[1mm]
   \hspace*{4cm}
           (A \equiv \pm \bpm 1&0\\0&1\epm  \mbox{ \em or } 
            A\equiv \pm \bpm0&-1\\1&0\epm \mbox{ \em mod } n)\}
          \end{array}         
\end{eqnarray*}
\end{thm}

\begin{proof} We proceed in two steps:\\

{\bf First step:}\\
We first determine $\stab(S)$.
We have:
\[ 
\begin{array}{l}
A \in \stab(S) \\
\hspace*{5mm} \ifff   A\cdot S = S \mbox{ with } S = \{b_1:=\bpm  p\\ q\epm,\; b_2:=\bpm -q\\p \epm,\; 
                                                 b_3:=\bpm -p\\-q\epm,\; b_4:=\bpm q\\-p \epm\}\\[2.5mm]
\hspace*{5mm} \ifff  \bar{A}\cdot S = S, 
            \mbox{ where $\bar{A}$ is the image of $A$ in $\slzwei(\ZZ/n\ZZ)$.} 
\end{array}
\]
Therefore $\stab(S)$ is a congruence group of level $n$.
We shall show that its image in $\slzwei(\ZZ/n\ZZ)$ is equal to:
\[\Gammaquer_S \; := \; \{\pm\bpm 1&0\\0&1\epm, \pm \bpm0&-1\\1&0\epm\}  
  \;\; \subseteq \;\;
  \slzwei(\ZZ/n\ZZ)\]
Observe that these four matrices  map $S$ to itself. Hence 
$\Gammaquer_S$ is contained in the image of $\stab(S)$.
Let now $\bar{A}$ be a matrix in the image, i.e. $\bar{A}\cdot S = S$. 
We have to show that $\bar{A}$ is in $\Gammaquer_S$.\\
By the definition of $\bar{B} \in \glzwei(\ZZ/n\ZZ)$ in (\ref{B}), 
we have $\bar{A}\cdot \bar{B} = (\bar{A}\cdot b_1, \bar{A}\cdot b_2)$, where
$(\bar{A}\cdot b_1, \bar{A}\cdot b_2)$ is the matrix whose first column is
$\bar{A}b_1$ and whose second column is $\bar{A}b_2$. 
Its determinant is equal to $\det(\bar{A})\cdot\det(\bar{B}) =
\det(\bar{B})$. 
Define $d := \det(\bar{B}) =  p^2 + q^2 \in (\ZZ/n\ZZ)^{\times}$.
We consider now the following four cases:\\[2mm]
{\em Case 1:}\;  $\bar{A}b_1 = b_1$.\\
$b_2$ is the only element of $S$ such that the determinant of the matrix 
$(b_1, b_2)$
is equal to $d$. {Hence it follows from the above} that 
$\bar{A}\cdot b_2 = b_2$.
Thus:\\[1mm]
$\bar{A}\cdot\bar{B} = \bar{B} \;\; 
 \stackrel{\bar{B} \in \mbox{\scs GL}_2(\ZZ/n\ZZ)}{\Longrightarrow} 
   \;\; \bar{A} = \bar{I}$, \; with $\bar{I}$ the
identity matrix in $\slzwei(\ZZ/n\ZZ)$.\\[2mm]
In the {\em other three cases}, namely $\bar{A}\cdot b_1 = b_2$,\; 
$\bar{A}\cdot b_1 = b_3$
and $\bar{A}\cdot b_1 = b_4$, we similarly obtain:
\[\bar{A} = \bpm 0&-1\\1&0\epm,\; \bar{A} = \bpm-1&0\\0&-1\epm \mbox{ and } 
  \bar{A} = \bpm 0&1\\-1&0\epm.\]
Hence, $\Gammaquer_S$ is the full image of $\stab(S)$ in $\slzwei(\ZZ/n\ZZ)$ and
$\stab(S)$ is the preimage of $\Gammaquer_S$ in $\slzwei(\ZZ)$.\\

{\bf Second step:}
{
For any $k \in \NN$ let $p_k: \slzwei(\ZZ) \to \slzwei(\ZZ/k\ZZ)$
be the natural projection.
We determine the Veech group $\Gamma(D_P)$ as subgroup of 
$\Gamma_S := p_n^{-1}(\Gammaquer_S)$.} \\[2mm]
From Corollary \ref{indexdpcor} and  Lemma \ref{two} it follows that 
\begin{equation}\label{chain}
  \Gamma(D_P)\;\; \subseteq\;\; \stab(S) \cap \Gamma_{u,u}\;\; = \;\;
  \Gamma_S \cap \Gamma_{u,u} \;\;=:\;\; \Gamma.
\end{equation}
We show that the index $[\Gamma_S:\Gamma]$ of $\Gamma$ in $\Gamma_S = \stab(S)$
is equal to 3. It then follows that we have equality in (\ref{chain}),
since $[\Gamma_S:\Gamma(D_P)]$ is also equal to 3 by Proposition \ref{indexprop}. 
This finishes the proof.\\
$\Gamma_S$ and $\Gamma$ are both congruence groups of level $2n$. Hence we may
as well show that their images in $\slzwei(\ZZ/2n\ZZ)$  differ by index 3, i.e.
that $[p_{2n}(\Gamma_S):p_{2n}(\Gamma)] = 3$.\\
Recall that  $\slzwei(\ZZ/2n\ZZ) \cong \slzwei(\ZZ/n\ZZ)\times \slzwei(\ZZ/2\ZZ)$,
since $\gcd(2,n) = 1$. Thus we have:
\begin{eqnarray*}
p_{2n}(\Gamma_S) 
  &=& \{A \in \slzwei(\ZZ/2n\ZZ)|\; A \equiv \pm I, \pm \bpm0&-1\\1&0\epm 
                                                      \mbox{ mod } n\}\\
  &\cong& \;\; p_n(\Gamma_S) \times  \slzwei(\ZZ/2\ZZ)
                    \;\; \subseteq \;\; \slzwei(\ZZ/n\ZZ) \times \slzwei(\ZZ/2\ZZ)\\
p_{2n}(\Gamma) &=& \{A \in \slzwei(\ZZ/2n\ZZ)|\;
    \begin{array}[t]{l}
    A \; \equiv \; \pm I, \pm \bpm0&-1\\1&0\epm \mbox{ mod } n \mbox{ and }\\[1mm]
    A \; \equiv \;I,  \bpm0&-1\\1&0\epm \mbox{ mod } 2 \}
    \end{array}\\
  &\cong&\;\; p_n(\Gamma_S) \times \{ I, \bpm0&-1\\1&0\epm\}
                     \;\; \subseteq \;\; \slzwei(\ZZ/n\ZZ) \times \slzwei(\ZZ/2\ZZ)
\end{eqnarray*}

Thus, $[\Gamma_S:\Gamma] = [p_{2n}(\Gamma_S):p_{2n}(\Gamma)) = 
       |\slzwei(\ZZ/2\ZZ)|/2 = 3$. 
\end{proof}

%% file: bild7q5.tex
\setlength{\unitlength}{0.7cm}
\begin{picture}(7,7)
\label{bildcase3}
\put(1,1){\framebox(1,1){}}
\put(1,2){\framebox(1,1){}}
\put(1,3){\framebox(1,1){}}
\put(1,4){\framebox(1,1){}}
\put(1,5){\framebox(1,1){}}
\put(2,1){\framebox(1,1){}}
\put(2,2){\framebox(1,1){}}
\put(2,3){\framebox(1,1){}}
\put(2,4){\framebox(1,1){}}
\put(2,5){\framebox(1,1){}}
\put(2,1){\framebox(1,1){}}
\put(2,2){\framebox(1,1){}}
\put(2,3){\framebox(1,1){}}
\put(2,4){\framebox(1,1){}}
\put(2,5){\framebox(1,1){}}
\put(3,1){\framebox(1,1){}}
\put(3,2){\framebox(1,1){}}
\put(3,3){\framebox(1,1){}}
\put(3,4){\framebox(1,1){}}
\put(3,5){\framebox(1,1){}}
\put(4,1){\framebox(1,1){}}
\put(4,2){\framebox(1,1){}}
\put(4,3){\framebox(1,1){}}
\put(4,4){\framebox(1,1){}}
\put(4,5){\framebox(1,1){}}
\put(5,1){\framebox(1,1){}}
\put(5,2){\framebox(1,1){}}
\put(5,3){\framebox(1,1){}}
\put(5,4){\framebox(1,1){}}
\put(5,5){\framebox(1,1){}}

\put(2.05,0.45){$P_0$}
\put(2.1,6.2){$P_0$}
\put(6.1,2.1){$P_1$}
\put(0.2,2.1){$P_1$}
\put(4.6,6.25){$P_2$}
\put(4.6,0.4){$P_2$}
\put(0.3,4.6){$P_3$}
\put(6.1,4.6){$P_3$}

\put(5.55,5.6){$B$}
\put(5.4,5.4){$\times$}
\put(2,1){\circle*{.3}}
\put(2,6){\circle*{.3}}
\put(6,2){\circle*{.3}}
\put(1,2){\circle*{.3}}
\put(5,6){\circle*{.3}}
\put(5,1){\circle*{.3}}
\put(1,5){\circle*{.3}}
\put(6,5){\circle*{.3}}

\put(0.2,.6){$O$}
\put(1,1){\circle*{.1}}
\put(1,1){\circle{.2}}

\put(3.6,3.6){$M$}
\put(3.5,3.5){\circle*{.1}}
\put(3.5,3.5){\circle{.2}}



\end{picture}